\documentclass{amsa}

\usepackage{amsmath,amssymb,amsfonts,amsthm}
\usepackage{mathtools}
\usepackage{fancyhdr,fancyvrb}
\usepackage{graphicx}
\usepackage{nextpage}
\usepackage{hyperref}
\usepackage{tikz}
\usepackage{tikz-3dplot}
\usepackage{mathrsfs}
\usepackage{bbm}
\usepackage{bm}
\usepackage{soul}

\allowdisplaybreaks

\newcommand{\compresslist}{
	\setlength{\itemsep}{1.5pt}
	\setlength{\parskip}{0pt}
	\setlength{\parsep}{0pt}
}

\newcommand{\R}[1]{\mathbb{R}^{#1}}

\newcommand{\bfx}{\textbf{x}}

\newcommand{\cC}{\mathcal C}

\newcommand{\de}{\mathrm d}

\newcommand{\loc}{\mathrm{loc}}

\newcommand{\rr}{\mathrm{r}}

\newcommand{\eps}{\varepsilon}

\newcommand{\abs}[1]{\left|#1\right|}

\renewcommand{\O}{\Omega}

\usepackage{soul}

\newcommand{\blu}[1]{\textcolor[rgb]{0,0,0}{#1}}

\renewcommand{\geq}{\geqslant}
\renewcommand{\leq}{\leqslant}
\newcommand{\wto}{\rightharpoonup}

\newcommand{\average}{{\mathchoice {\kern1ex\vcenter{\hrule
height.4pt width 8pt depth0pt}
\kern-11pt} {\kern1ex\vcenter{\hrule height.4pt width 4.3pt
depth0pt} \kern-7pt} {} {} }}

\mathchardef\emptyset="001F

\newtheorem{defin}{Definition}[section]
\newtheorem{remark}[defin]{Remark}
\newtheorem{theorem}[defin]{Theorem}

\newcommand{\editor}{Xxxxxxxx, XXXX}
\newcommand{\rcptDate}{Xxxx XX, XXXX}

\begin{document}
\label{page:t}
\thispagestyle{plain}
\title{GEOMETRICALLY CONSTRAINED WALLS \\ 
IN THREE DIMENSIONS
}
\author{Riccardo Cristoferi, Gabriele Fissore}
\affiliation{Department of Mathematics - IMAPP, Radboud University, Nijmegen, The Netherlands}
\email{\{riccardo.cristoferi,gabriele.fissore\}@ru.nl}
\sauthor{Marco Morandotti}
\saffiliation{Dipartimento di Scienze Matematiche ``G.~L.~Lagrange'', Politecnico di Torino, Italy}
\semail{marco.morandotti@polito.it}
%
\footcomment{
GF was partially supported under NWO-OCENW.M.21.336. This study was carried out within the \emph{Geometric-Analytic Methods for PDEs and Applications} project (2022SLTHCE, CUP E53D23005880006), funded by European Union -- Next Generation EU within the PRIN 2022 program (D.D. 104 - 02/02/2022 Ministero dell’Universit\`{a} e della Ricerca). This manuscript reflects only the authors’ views and opinions and the Ministry cannot be considered responsible for them.
}
\maketitle
\noindent
{\bf Abstract.} \indent
We study geometrically constrained magnetic walls in a three dimensional geometry where two bulks are connected by a thin neck.
Without imposing any symmetry assumption on the domain, we investigate the scaling of the energy as the size of the neck vanishes.
We identify five significant scaling regimes, for all of which we characterise the energy scaling and identify the asymptotic behaviour of the domain wall.
Finally, we notice the emergence of sub-regimes that are not present in the previous works due to restrictive symmetry assumptions.
\newpage

\section{Introduction} \label{sec:intro}
\indent 

\begin{figure}[h]
\begin{center}
\includegraphics[scale=.3]{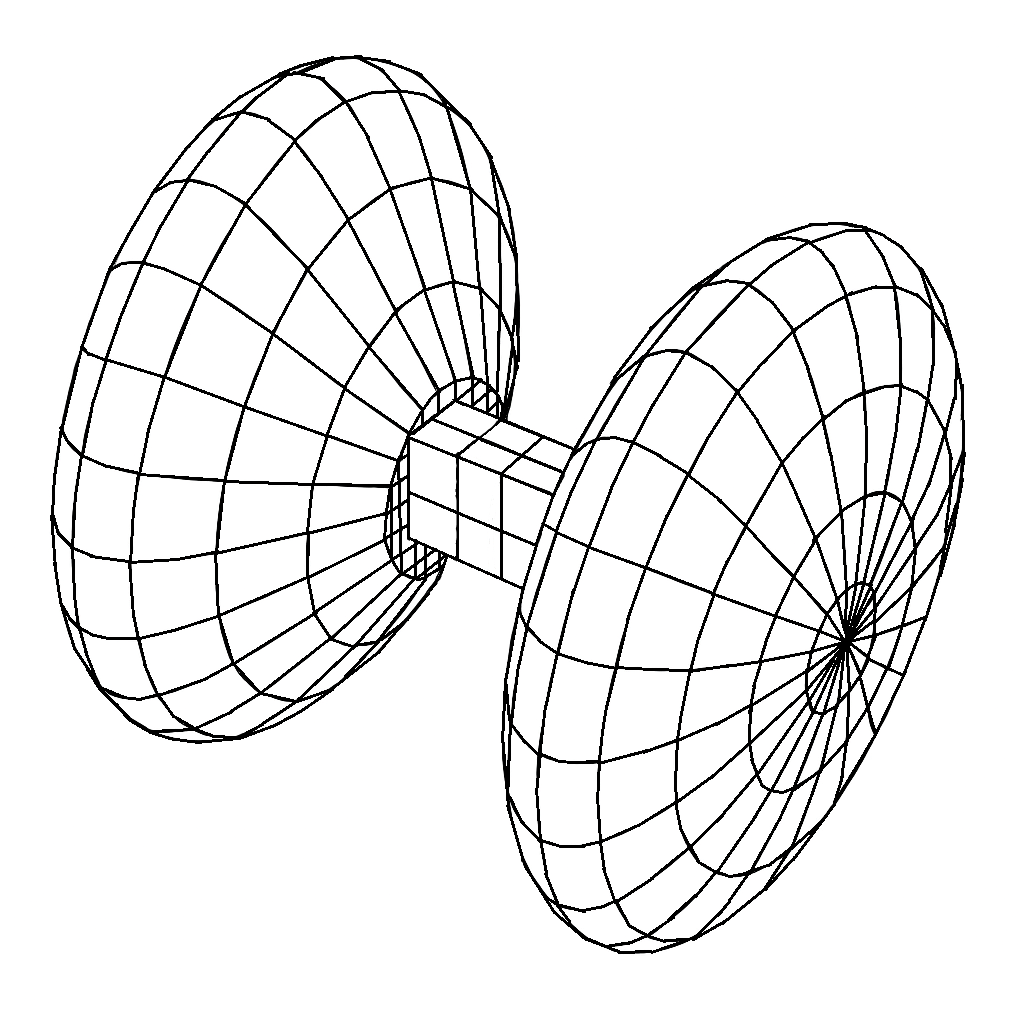}
\caption{A pictorial representation of a typical domain of interest.}
\label{figone}
\end{center}
\end{figure}

\emph{Magnetic domain walls} are regions in which the magnetisation of a material changes from one value to another one. 
In the presence of extreme geometries, such as that of a dumbbell-shaped domain (see Figure~\ref{figone}), the magnetic wall is more likely to be found in or around the neck; in this and similar geometry-driven situations, one usually speaks of \emph{geometrically constrained walls}, to stress the fact that the domain shape plays a pivotal role in the localisation of the transition region of the magnetisation, for instance, when prescribing it in the bulky parts of the dumbbell.

The study of the behaviour of the magnetisation in a dumbbell-shaped domain is relevant in micro- and nano-electronics application, where the neck of the dumbbell models magnetic point contacts.
We refer the reader to \cite{CheYan10,JubAllBis04,Klaui_2008,Mol_etal02,Sasaki_etal12,Tatara_etal99} for an incomplete list of applications and contexts of relevance of geometrically constrained walls.
If one imposes two different values of the magnetisation, one in each of the two macroscopic components, a transition is expected in the vicinity of the neck, as initially observed by Bruno in \cite{Bru99}: if the neck is small enough, so that the geometry of the material varies drastically when passing from one bulk to the other, it can play a crucial role in determining the location of the magnetic wall, by influencing the mere minimisation of the magnetic energy. 
Three scenarios are to be considered: the transition may happen either completely inside the neck, or partly inside and partly outside the neck, or completely outside the neck. 

The model features a sufficiently smooth potential which is minimal at the imposed values of the magnetisation in the bulk parts of the dumbbell, and a gradient term penalising transitions; the two are competing as soon as the values of the magnetisation in the bulks are not the same. 
To model the geometries that are of interest in the applications, we will consider an infinitesimally small neck, whose size is determined by three parameters $\eps, \delta, \eta>0$:
\begin{equation}\label{eq_neck_ep}
N_{\eps}\coloneqq\{\bfx=(x,y,z)\in\R{3}: |x|\leq \eps,|y|<\delta,|z|<\eta\},
\end{equation}
with the understanding that all three of them vanish when $\eps\to0$, that is $\delta=\delta(\eps)\to0$ and $\eta=\eta(\eps)\to0$, as $\eps\to0$. 
The full domain is described by
\begin{equation}\label{eq_domain_eps}
\Omega_{\eps} \coloneqq \Omega_{\eps}^{\ell} \cup N_{\eps} \cup \Omega_{\eps}^{\rr},
\end{equation}
where $\Omega_{\eps}^{\ell}=\Omega^{\ell}-(\eps,0,0)$ and $\Omega_{\eps}^{\rr}=\Omega^{\rr}+(\eps,0,0)$, for certain open sets $\Omega^{\ell}\subset\{x<0\}$ and $\Omega^{\rr}\subset\{x>0\}$ such that $0\in\partial\Omega^{\ell}\cap\partial\Omega^{\rr}$.
This geometry makes the $x$ direction the preferred one, whereas the $y$- and $z$-direction can be interchanged upon a change of coordinates; this motivates the fact that we will use, throughout the work, the subscript $\eps$ alone as an indication of the smallness of the neck.

We are interested in understanding the asymptotic behaviour, as $\eps\to0$, of stable critical points (see Definition \ref{def:loc_min}) of the energy
\begin{equation}\label{004}
	F(u,\Omega_\eps) \coloneqq 
	\frac{1}{2} \int_{\Omega_\eps} \abs{\nabla u(\textbf{x})}^2\,\de \textbf{x} + \int_{\Omega_\eps} W(u(\textbf{x}))\,\de \textbf{x}\,, 
\end{equation}
defined for $u\in H^1(\Omega_\eps)$, where $\de \textbf{x} = \de x\de y \de z$ and 
\begin{equation}\label{eq_hp_W}
\begin{cases}
\text{$W\colon\R{}\to[0,+\infty)$ is a function of class $\cC^2$ such that}\\[5pt]
\text{$W^{-1}(0)=\{\alpha,\beta\}$ for some $\alpha<\beta$ and $\displaystyle \lim_{|t|\to+\infty} W(t)=+\infty$.}
\end{cases}
\end{equation}
In \eqref{004}, the function $u$ represents a suitable quantity related to the magnetisation field defined on $\Omega_\eps$ and the potential $W$ favours the values $u(\textbf{x})=\alpha$ and $u(\textbf{x})=\beta$ for the magnetisation, corresponding to those to be imposed in the bulks.
Here, the competition between the potential and the gradient terms is significantly influenced by the geometry of $\Omega_\eps$.
The energy \eqref{004} was considered in \cite{Bru99} as a simplified model for studying the magnetisation inside a thin dumbbell domain under the assumption that the magnetic field is of the form
\begin{equation}\label{eq:Bruno}
\bold{m}(\bold{x})=\left(0, \cos(u(x)), \sin(u(x)) \right).
\end{equation}
Despite this simplifying assumption, the mathematical analysis is rich enough to exhibit non-trivial behaviours of the magnetisation.

\subsection{Related literature}
The body of literature on this problem counts many physical contributions stemming from Bruno's work \cite{Bru99} and a few mathematical items, which we are going to briefly review to give a context to our novel results. 

In \cite{Bru99}, Bruno considers the special geometry of a thin ($0<h\ll1$) three-dimensional wall $\Omega=S\times(-h,h)$, where $S$ is a planar region with a dumbbell shape whose neck is located at the origin and the bulks are in $\{x<0\}$ and $\{x>0\}$.
He assumes that the preferred directions of the magnetisation vector are $\mathbf{m}=(0,0,\pm1)$ and makes the Ansatz \eqref{eq:Bruno} that it varies only in the $y$-$z$ plane, as a function of the $x$-coordinate alone.
The energy that Bruno minimises is the one usually describing Bloch walls and turns out to be a functional of $u$ alone, with two emerging length scales when imposing that $\mathbf{m}\approx(0,0,\pm1)$ in $\{\pm x>0\}$: one driven by the shape $S$ of the domain, the other one dictated by the physical parameters entering the expression of the energy.
Despite Bruno's insightful intuition, the special form of the magnetisation has some limitations, some of which were removed (for instance, by allowing $\mathbf{m}$ to vary also in the $x$-$z$ plane and considering fully three dimensional geometries) in \cite{MolOsiPon02}.

Among the mathematical literature on this topic, we mention \cite{ArrCar04, ArrCarLC06, ArrCarLC09a, ArrCarLC09b, CDLLMur08, HalVeg84, Jim84, Jim88, Jim04} as far as the PDE aspect is concerned, and \cite{CabFreMorPer01, CDLLMur04, CDLLMur04b,  RubSchSte04} for variational approaches. 
Finally, we discuss two more recent contributions in detail, since the results we obtain are related to them.

In the work by Kohn \& Slastikov \cite{KohSla}, the problem is studied in the full three-dimensional setting, with the assumption that the geometry be axisymmetric: the dumbbell $\Omega_{\epsilon}$ is a rotation body around the $x$ axis, so that the shape parameters of the neck are essentially its length $\eps$ and its radius $\delta$.
\blu{In their work, they also consider an $\varepsilon$-dependent parameter $\gamma_\varepsilon$ in front of the potential term, and study the behavior of stable critical points also in this case. Moreover, they consider also more general geometries for the neck $N_\varepsilon$.}
By taking advantage of a scale-invariant Poincar\'{e} inequality for Sobolev functions and by reducing the problem to a one-dimensional variational one, the authors proved the existence of three possible regimes, according to the value of the limit $\lim_{\eps\to0} \delta/\eps=\lambda\in[0,+\infty]$ and singled out a \emph{thin neck} regime ($\lambda=0$), a \emph{normal neck} regime ($\lambda\in(0,+\infty)$), and a \emph{thick neck} regime ($\lambda=+\infty$).
In the first case the transition happens entirely inside the neck and is an affine function of the $x$-coordinate, in the second case the transition happens across the neck, partially inside and partially outside, depending on the value of~$\lambda$, whereas in the third case the transition happens entirely outside the neck.
These behaviours are found by studying the energy of particular competitors (essentially, an affine transition inside the neck and a harmonic transition in a spherical shell just outside the neck) and then rescaling the minimiser in the vicinity of the neck.

In the works by Morini \& Slastikov \cite{MorSla12, MorSla15} the same problem was addressed in the case of magnetic thin film, that is when the domain has the shape of a dumbbell, but it is two-dimensional, that is, in Bruno's setting in the limit as $h\to0$.
Mathematically speaking, the endeavour is more difficult on two accounts: the scale-invariant Poincar\'{e} inequality is not available in dimension two, and the problem loses its variational character.
Methods that are typical from the study of PDE's were employed to construct suitable barriers to estimate the solutions. 
Moreover, due to the slow decay of the logarithm (the fundamental solution to Laplace's equation in two dimensions), sub-regimes became available in addition to the thin, \blu{normal}, and thick neck regimes already analyzed by Kohn \& Slastikov: the sub-critical, critical, and super-critical thin neck regimes were found according to the value of the limit $\lim_{\eps\to0} (\delta|\ln\delta|)/\eps\in[0,+\infty]$, displaying a richer variety.
In the case of sub-regimes, the rescaling of the minimisers to study their asymptotic behaviour is not trivial; nonetheless, the authors managed to characterise the profiles as the unique solutions to certain PDE's where the boundary conditions track the expected asymptotic behaviour.

Both in Kohn \& Slastikov's and in Morini \& Slastikov's papers the technique involves two steps: (i) estimate the energy of minimisers to understand if the wall is located all inside, all outside, or across the neck, and
(ii) rescaling the whole domain $\Omega_\eps$ to an appropriate $\Omega_\infty$ in a way that either a variational problem or a PDE can be studied in $\Omega_\infty$, which brings to the attention that the boundary $\partial\Omega_\infty$ must be a set in which boundary conditions can be prescribed.


\section{Main result}

\subsection{Set up of the problem}\label{sec:results}
We study a mathematical model to characterise magnetic domain walls in a three-dimensional dumbbell-shaped domain (see Figure~\ref{figone}). 
The two bulks are modelled by two bounded, connected, open \blu{Lipschitz} sets $\Omega^\ell,\Omega^\rr\subset\mathbb{R}^3$ such that 
\begin{itemize}\compresslist
	\item[(H1)] the origin $(0,0,0)$ belongs to $\partial\Omega^\ell \cap \partial\Omega^\rr$;
	\item[(H2)] $\Omega^\ell \subset \{x<0\}$ and $\Omega^\rr \subset \{x>0\}$;
	\item[(H3)] there exists $r_0>0$ such that $(\partial\Omega^\ell) \cap B_{r_0}(0,0,0)$ and $(\partial\Omega^\rr) \cap B_{r_0}(0,0,0)$ are contained in the plane $\{x=0\}$, i.e., the bulks are flat and vertical near the origin, where the conjunction with the neck will be located.\footnote{We point out that this assumption is made for mere convenience and it does not affect the generality of our results. While allowing the reader to focus on the main qualitative geometrical properties of the domain, it can be removed following the strategy outlined in \cite{MorSla12}.}
\end{itemize}
We let $\eps>0$ and define the neck region as in \eqref{eq_neck_ep}, so that the dumbbell-shaped domain $\Omega_\eps$ is defined as in \eqref{eq_domain_eps}, where
$\Omega_\eps^\ell = \Omega^\ell - (\eps,0,0)$ and $\Omega_\eps^\rr = \Omega^\rr + (\eps,0,0)$,
We notice that $\Omega_\eps$ is a bounded, connected, open set with Lipschitz boundary.

We now give the relevant definitions of critical points and isolated local minimiser for the functional $F(\cdot,\Omega_\eps)$ introduced in \eqref{004}.
\begin{defin}\label{def:crit_point}
We say that a function $u_\eps\in H^1(\O_\eps)$ is a \emph{critical point of $F(\cdot,\Omega_\eps)$} if it is a weak solution to the system
\begin{align*}
	\begin{cases}
		\Delta u_\eps= W'(u_\eps)\quad&\text{in}\ \O_\eps,\\[5pt]
		\displaystyle\frac{\partial u_\eps}{\partial \nu}=0\quad&\text{on}\ \partial\O_\eps.
	\end{cases}
\end{align*}
\end{defin}

\begin{defin}\label{def:loc_min}
For $\eps>0$, let $u_\eps\in H^1(\Omega_\eps)$ be a critical point of $F(\cdot,\Omega_\eps)$.
We say that the family $(u_\eps)_\eps$ is an \emph{admissible family of nearly locally constant critical points} if
\begin{enumerate}\compresslist
\item[(a)] there exists $\bar\eps>0$ such that $\sup \big\{ \lVert u_\eps \rVert_{\infty}: \eps\in(0,\bar\eps] \big\} \eqqcolon \overline{M}<+\infty$;
\item[(b)] $\lVert u_\eps-\alpha\rVert_{L^1(\Omega_\eps^\ell)}\to0$ and     $\lVert u_\eps-\beta\rVert_{L^1(\Omega_\eps^\rr)}\to0$, as $\eps\to0$.
\end{enumerate}
Moreover, we say that $(u_\eps)_\eps$ is an \emph{admissible family of local minimisers} if it satisfies, additionally, 
\begin{enumerate}\compresslist
\item[(c)] there exist $\eps_0>0$ and $\theta_0>0$ such that for $\eps\in(0,\eps_0]$ we have
$$\text{$F(v,\Omega_\eps)\geq F(u_\eps,\Omega_\eps)$ \quad for all $v\in H^1(\Omega_\eps)$ such that $0<\lVert v-u_\eps\rVert_{L^1(\Omega_\eps)}\leq\theta_0$\,.}$$
\end{enumerate}
\end{defin}

\blu{
\begin{remark}
In \cite[Definition 3.1]{MorSla15} the authors introduced the same notion by requiring the $L^2$ convergence in (b). Due to (a), our definition is equivalent to theirs. We chose to use the $L^1$ convergence because it is a weaker notion to be checked.
\end{remark}
}

Regarding the existence of minimisers, \cite[Theorem~3.1]{KohSla} can easily be adapted to our setting.
\begin{theorem}
For $\eps>0$, let $u_{0,\eps}\colon\Omega_\varepsilon\to\R{}$ be defined as
\begin{align*}
	u_{0,\eps}(\bold{x})\coloneqq\begin{cases}
		\alpha & \text{if } \bold{x}\in\Omega^\ell_\eps\,,\\[5pt]
		\displaystyle \frac{\alpha+\beta}{2}  & \text{if } \bold{x}\in  N_\eps\,,\\[5pt]
		\beta & \text{if } \bold{x}\in\Omega^\rr_\eps\,.
	\end{cases}
\end{align*}
If $u_\eps\in H^1(\Omega_\eps)$ is such that $F(u_\eps,\Omega_\eps)\leq F(v,\Omega_\eps)$ for every $v\in B_\eps$, where
\begin{equation}\label{Beps}
	B_\eps\coloneqq\{u\in H^1(\O_\eps):||u-u_{0,\eps}||_{L^2(\O_\eps)}\leq d,\, \|u\|_{L^\infty(\Omega_\eps)}<\infty  \},
\end{equation}
with $d<\min\{\blu{\abs{\alpha}}|\O^\ell|^{1/2},\blu{\abs{\beta}}|\O^\rr|^{1/2}\}$, then the family $(u_\eps)_\eps$ is an admissible family of local minimisers according to Definition~\ref{def:loc_min}, and $\lVert u_\eps - u_{0,\eps} \rVert_{L^2(\Omega_\eps)}\to0$, as $\eps\to0$.
\end{theorem}
Unlike \cite{KohSla}, we do not assume axial symmetry of the domain and this results in a richer variety of regimes.
In particular, we find that some of these regimes admit sub-regimes, as was discovered for magnetic thin films in \cite{MorSla12, MorSla15}.
We discuss all the possible cases in the next section.


\subsection{Heuristics}\label{sec:heuristics}

In this section, we show how to heuristically guess where the main part of the energy concentrates, just by looking at the asymptotic relationships between the three geometric parameters $\eps,\delta,\eta$.

First of all, we note that, given the privileged role of the parameter~$\eps$, it is trivial to see that the roles of~$\delta$ and~$\eta$ can be interchanged upon switching the coordinate axes~$y$ and~$z$. 
The regimes investigated in \cite{KohSla} corresponds to the cases where $\delta\sim\eta$.
Therefore, here we limit ourselves to consider the other following regimes:
\begin{enumerate}
	\item[$(i)$] \textit{Super thin}: $\varepsilon\gg\delta\gg\eta$;
	\item[$(ii)$] \textit{Flat thin}: $\eps\approx\delta\gg\eta$;
	\item[$(iii)$] \textit{Window thick}: $\delta \gg\eta\gg \eps$;
	\item[$(iv)$] \textit{Narrow thick}: $\delta\gg\eps\approx\eta$;
	\item[$(v)$] \textit{Letter-box}: $\delta\gg\eps\gg\eta$.
\end{enumerate}
We now want to guess where the transition will happen: completely inside, completely outside, or in both regions.
To understand this, we reason as follows.
First of all, we expect the main part of the energy to be the Dirichlet integral.
Therefore, we consider two harmonic functions that play the role of competitors for the minimisation problem 
$$\min \{ F(v,\Omega_\eps)\,: v\in B_\eps\};$$
one where the transition from $\alpha$ to $\beta$ happens inside the neck, and the other one where it happens only outside (and the competitor is constant inside the neck).
We then compare their energies (whose computations will be carried out in Section \ref{sec:competitors}) to get a guess of where the transition will occur.
The first harmonic function will be referred to as the \emph{affine competitor}, and has energy of order
\[
\text{Energy of the affine competitor } = \frac{\delta\eta}{\varepsilon}.
\]
The second harmonic function will be referred to as the \emph{elliptical competitor}, and has energy of order
\[
\text{Energy of the elliptical competitor } = \frac{\delta}{\abs{\ln(\eta/\delta)}}.
\]
When one of the two energies is dominant with respect to the other, it is clear where we expect the transition to happen.
In the case they are of the same order, we guess that the transition is both inside and outside the neck. This will be later confirmed by rigorous analysis.

The comparison of the energies of the two harmonic competitors leads to the following heuristics:
\begin{enumerate}
	\item[$(i)$] \textit{Super thin neck}: In this case, we have 
	\begin{align*}
		\frac{\delta\eta}{\eps} \frac{|\ln(\eta/\delta)|}{\delta}=\frac{\delta}{\eps}\Big(\frac{\eta}{\delta} \abs{\ln(\eta/\delta)}\Big)\to 0,
	\end{align*}
	as $\eps\to0$. Then we expect the transition to happen inside $N_\eps$.
	\item[$(ii)$] \textit{Flat thin neck}: In this case, we obtain
	\begin{align*}
		\frac{\delta\eta}{\eps} \frac{|\ln(\eta/\delta)|}{\delta}=\frac{\eta}{\delta}\abs{\ln(\eta/\delta)}\to 0,
	\end{align*}
	thus the transition is occurring inside $N_\eps$.
	\item[$(iii)$]\textit{Window thick neck}: The comparison of the energies of the harmonic competitors gives
	\begin{align*}
		\frac{\delta\eta}{\eps} \frac{\abs{\ln(\eta/\delta)}}{\delta}=\frac{\eta}{\eps}	\abs{\ln(\eta/\delta)}\to+\infty,
	\end{align*}
as $\varepsilon\to0$.
	The transition is expected to happen entirely outside the neck.
	\item[$(iv)$] 	
	 \textit{Narrow thick neck}: In this case, we have
	\begin{align*}
		\frac{\delta\eta}{\eps} \frac{\abs{\ln(\eta/\delta)}}{\delta}=\abs{\ln(\eta/\delta)}\to+\infty,
	\end{align*}
	as $\eps\to0$. Therefore, we expect the transition to happen outside the neck.
	\item[$(v)$]
 \textit{Letter-box neck}: In this case, we have the presence of sub-regimes. Indeed, the comparison of the orders of the energies of the harmonic competitors gives
	\begin{align*}
		\frac{\delta\eta}{\eps} \frac{\abs{\ln(\eta/\delta)}}{\delta}=\frac{\eta}{\eps}\abs{\ln(\eta/\delta)},
	\end{align*}
whose asymptotic behaviour is not clear.
Therefore, we have to consider the following sub-regimes:
	\begin{enumerate}
		\item[$(a)$] \textit{Sub-critical letter-box neck}, when
		\begin{align*}
			\frac{\delta\eta}{\eps} \frac{\abs{\ln(\eta/\delta)}}{\delta}\to 0,
		\end{align*}
		as $\eps\to0$. In such a case, we expect the transition to happen inside the neck.
		\item[$(b)$] \textit{Critical letter-box neck}, when
		\begin{align*}
			\frac{\delta\eta}{\eps} \frac{\abs{\ln(\eta/\delta)}}{\delta}\to \ell \in(0,+\infty),
		\end{align*}
		as $\eps\to 0$.
In this case, we expect the transition to happen both inside and outside the neck.
		\item[$(c)$] \textit{Super critical letter-box neck}, when
		\begin{align*}
			\frac{\delta\eta}{\eps} \frac{\abs{\ln(\eta/\delta)}}{\delta}\to +\infty,
		\end{align*}
		as $\eps\to0$. Here, the transition is expected to happen outside the neck.
	\end{enumerate}	
\end{enumerate}

\subsection{Novelty of the present contribution and main results}
In this paper, we study the full three-dimensional case of the problem with no symmetry assumption on the shape of the neck: it will be a  parallelepiped as in \eqref{eq_neck_ep} with all three dimensions independent from one another and all vanishing to zero as $\eps\to0$.
When considering the mutual rate of convergence to zero of the three parameters $\eps,\delta,\eta$, we single out five regimes that do not emerge in the analysis of Kohn \& Slastikov, and for each of them we study the energy scaling. 
We notice that in our three-dimensional setting the scale-invariant Poincar\'{e} inequality is not always available.
This inequality ensures that, given a \blu{bounded, connected} open set $A\subset\R{3}$, there exists a constant $C>0$ such that
\[
\left(\int_{\lambda A} \left|u\left(\frac{\textbf{x}}{\lambda}\right)
	- \overline{u}_{A}\right|^6\ \de \textbf{x} \right)^{\frac{1}{6}}
\leq
C \left( \int_{\lambda A} \left|\nabla u\left(\frac{\textbf{x}}{\lambda}\right)\right|^2\ \de \textbf{x} \right)^{\frac{1}{2}},
\]
for all $\lambda>0$, and all $u\in H^1(A)$.
Here, $\overline{u}_{A}$ denotes the average of $u$ in $A$.
Note that the argument to get such inequality is the same to guess the conjugate exponent in the Gagliardo-Nierenberg inequality.
Despite that, we are able to investigate the behaviour of local minimisers and the associated rescaled limiting problem, which possesses a variational structure in every regime.

Due to the peculiar geometry of our problem induced by the mismatch between $\eta$ and $\delta$, namely $\eta\ll\delta$, the cross section of the junction of the neck with the bulks is a rectangle with a very large aspect ratio; this allows us to find an ellipsoidal competitor carrying less energy than the spherical one proposed in the previous works.
As it depends on $\abs{\ln(\eta/\delta)}$, the energy scaling turns out to exhibit sub-regimes in some of these cases, as described in Section~\ref{sec:heuristics}.

The main achievements of the paper are the following.
\begin{itemize}
\item[(A1)] For all of the above-mentioned regimes, we identify where the transition happens.
More precisely, we find sequences $(\varrho_\varepsilon)_\eps$ depending explicitly on the parameters $\eps, \delta, \eta$, with $\varrho_\eps\to+\infty$, as $\eps\to 0$, such that
\begin{align*}
\lim_{\eps\to0} \varrho_\eps F_\eps(u_\eps,\O_\eps) &=: \kappa\in (0,+\infty),\\[5pt]
\lim_{\eps\to 0} \varrho_\eps F_\eps(u_\eps,N_\eps) &=: \kappa_N\in [0,\kappa].
\end{align*}
Their interpretation is the following: $\kappa$ is the asymptotic energy in the whole domain, and $\kappa_N$ that in the neck.
Therefore, if $\kappa_N=\kappa$, we say that the transition happens entirely inside the neck, if $\kappa_N=0$, the transition happens entirely outside the neck, while if $\kappa_N\in(0,\kappa)$, the transition happens both inside and outside the neck.
In particular, we rigorously justify the expectations derived from the above heuristics.
\item[(A2)] In all of the regimes, we consider the profile in the region where the transition happens.
We identify a proper rescaling of the independent variables that allow us to prove that such rescaled profile converges to a solution of a Dirichlet energy in a limiting domain with suitable boundary conditions.
Only in the critical letter-box regime, we need to assume convergence to a limiting profile (see Theorem \ref{Theoremcriticalletterbox}(i)), and we prove that the latter solves a minimisation problem (see Remark \ref{rem:critical}).
In all cases, local minimisers will converge to a constant in the region where the transition does not happen.
\end{itemize}

We refer the reader to Section \ref{sec:analysis} for the precise statements and proofs of these results.

\begin{remark}\label{rem:critical}
The reason why in the critical letter-box regime we cannot prove compactness of a sequence of local minimisers, is the following. We do expect to see part of the transition inside the neck. 
Therefore, we rescale the local minimiser $u_\eps$ as $v_\eps(x,y,z)\coloneqq u_\eps(\eps x, \delta y, \eta z)$.
In such a way, we get that
\[
\frac{\eps}{2\delta\eta}\int_{N_\eps}\abs{\nabla u_\eps}^2\ \de\bold{x}
= \frac{1}{2}\int_{[-1,1]^3}\Big((\partial_x v_\eps)^2+\frac{\eps^2}{\delta^2}(\partial_y v_\eps)^2+\frac{\eps^2}{\eta^2}(\partial_z v_\eps)^2\Big)\ \de\bold{x}.
\]
The left-hand side is bounded thanks to the energy of the affine competitor.
Unfortunately, since in this regime $\eps\ll\delta$, we do not get a lower bound of the $y$-derivative of the function $v_\eps$, even if we prove that each limit of the sequence $(v_\eps)_\eps$ will only depend on the first variable.
\end{remark}

\section{Competitors}\label{sec:competitors}
\indent

The goal of this section is to build two harmonic competitors and to compute the order of their energies. For clarity, we present the affine and elliptic competitors separately. However, at the end of the section, they are mixed together in a more general way.

\subsection{Affine competitor}
Here we build the affine competitor inside the neck and we compute its energy.
Let $A,B\in\R{}$, and define the affine function $\xi_\varepsilon \colon \R{3}\to\R{}$ as
\begin{align}\label{affinecompetitor}
	\xi_\varepsilon(\bfx)\coloneqq\begin{cases}
		A\quad & \text{if $\bold{x}\in\Omega^\ell_\eps$\,,}\\[3pt]
			\displaystyle\frac{B-A}{2\eps}x+\frac{B+A}{2}\quad & \text{if $\bold{x}\in N_\eps$\,,}\\[5pt]
		B\quad& \text{if $\bold{x}\in\Omega^\rr_\eps$\,.}
	\end{cases}
\end{align}
Then, we have that
\begin{align}\label{AffineBound}
 \frac{1}{2}\int_{N_\varepsilon} |\nabla \xi_\varepsilon|^2 \,\de \mathbf{x}
 =\frac{\delta\eta}{\eps}(B-A)^2.
\end{align}

\subsection{Prolate competitor}\label{sec:elliptic_competitor}
In \cite{KohSla}, the authors built a harmonic competitor by imposing boundary condition on \emph{half-spheres} centred at the edges of the neck. The choice of the spherical geometry was dictated by the fact that the authors required $\delta=\eta$.
In our case, the geometry will be that of an ellipsoid, suggested by the fact that one of the parameters $\delta$ and $\eta$ is larger than the other.

In order to define our prolate competitor, we first need to introduce the so-called \emph{prolate spheroidal coordinates}.
Consider, for $a>0$, the change of coordinates $(x,y,z) =\blu{ \Psi_a}(\mu, \nu, \varphi)$, given by
\begin{equation}\label{EllipsHConstr}
	\begin{cases}
		x = a \sinh \mu \sin \nu \cos \varphi, \\
		y =  a \cosh \mu \cos \nu, \\
		z =  a \sinh \mu \sin \nu \sin \varphi, \\
	\end{cases}
\end{equation}
where $(0, \pm a, 0)$ are the coordinates of the foci, $\varphi\in [0,2\pi]$ is the polar angle, $\nu\in[0,\pi]$ is the azimuthal angle, and $\mu>0$.
The Jacobian matrix and the Jacobian determinant of the transformation $\Psi_a$ are, respectively,
\begin{equation}\label{jacobian}
J_{\Psi_a}(\mu,\nu,\varphi)=
\begin{pmatrix}
a\cosh\mu \sin\nu \cos\varphi & a\sinh\mu \cos\nu \cos\varphi & - a\sinh\mu \sin\nu \sin\varphi \\
a\sinh\mu \cos\nu & -a\cosh\mu \sin\nu & 0\\
a\cosh\mu \sin\nu \sin\varphi & a\sinh\mu \cos\nu \sin\varphi & a \sinh\mu \sin\nu \cos\varphi
\end{pmatrix}
\end{equation}
and
\begin{equation}\label{jacobian_determinant}
\det(J_{\Psi_a}(\mu,\nu,\varphi))=-a^3 \sinh\mu \sin\nu (\sin^2\nu+ \sinh^2\mu).
\end{equation}
For $M>0$, define the ellipsoid
\begin{equation}\label{ellipsoid}
E(a,M)\coloneqq\{\Psi_a(\mu,\nu,\varphi): \mu<2M\}.
\end{equation}
Moreover, we need consider the left and the right halves of the set $E_M$ translated at the edges of the neck. Namely, we consider the open sets
\[
	E_\eps^\ell(a,M) \coloneqq (E(a,M)\cap\{x<0\}) - (\eps, 0, 0),
\]
and
\[
	E_\eps^\rr(a,M) \coloneqq (E(a,M)\cap\{x>0\}) + (\eps, 0, 0).
\]
Note that if $\blu{a\cosh(2M)<r_0}$, where $r_0>0$ is given by assumption (H3), we get that $E_\eps^\ell(a,M)\subset \Omega_\eps^\ell$, and that $E_\eps^\rr(a,M)\subset \Omega_\eps^\rr$.
For $0<m<M$, we define the function $\xi_\eps\colon\R{3}\to\R{}$ as
\begin{equation}\label{testfootball}
	\xi_\eps(\mathbf{x})=\xi_\eps(x,y,z)\coloneqq
	\begin{cases}
		\alpha
			& \text{ in }
				\Omega_\eps^\ell\setminus \overline{E}_\eps^\ell(a,M)\,, \\[5pt]
		\displaystyle \frac{\alpha+\beta}{2}-h(x+\eps,y,z)
			& \text{ in } E_\eps^\ell(a,M)\,, \\[5pt]
		\displaystyle \frac{\alpha+\beta}{2}
			& \text{ in } N_\eps\,, \\[5pt]
		\displaystyle \frac{\alpha+\beta}{2}+h(x-\eps,y,z)
			& \text{ in } E_\eps^\rr(a,M)\,, \\[5pt]
		\beta & \text{ in } \Omega_\eps^\rr\setminus \overline{E}_\eps^\rr(a,M)\,,
	\end{cases}
\end{equation}
where $h \colon E(a,M) \setminus E(a,m)\to\R{}$ is the solution to the problem
\begin{equation*}
	\begin{cases}
		\Delta h=0 & \text{in } E(a,M) \setminus \overline{E}(a,m)\,, \\[5pt]
		\displaystyle h = \frac{\beta-\alpha}{2} & \text{on } \partial E(a,M)\,, \\[5pt]
		h = 0 & \text{on } \partial E(a,m)\,.
	\end{cases}
\end{equation*}
Now, our goal is to find the function $h$ explicitly and to estimate, asymptotically, its Dirichlet energy.
We look for the solution in the form $\blu{h(x,y,z)= \tilde h(\mu)}$.  Then the Laplacian in prolate spheroidal coordinates is given by
\begin{equation*}
	\Delta \tilde{h}(\mu) = \frac{1}{a^2 (\sin^2 \nu + \sinh^2 \mu)} \left(\tilde{h}_{\mu \mu}(\mu) + (\coth \mu) \tilde{h}_\mu(\mu) \right) = 0,
\end{equation*}
or equivalently	$(\sinh \mu \tilde{h}_\mu)_\mu = 0$.
It follows that
\begin{equation}\label{larichiamiamo}
	\tilde{h}_\mu(\mu) = \frac{c}{\sinh \mu},
\end{equation}
and thus
\begin{equation}\label{hmu}
	\tilde{h}(\mu) = c \ln |k \tanh(\mu/2)|.
\end{equation}
Enforcing the boundary conditions
\begin{equation*}
	\tilde{h}(2M) = \frac{\beta-\alpha}{2}\quad\text{and}\quad \tilde{h}(2m) = 0
\end{equation*}
yields
\begin{equation}\label{eq:BS_elliptic_competitor}
	k = \frac1{\tanh m} \qquad\text{and}\qquad
	c = \frac{\beta-\alpha}{2\ln\Big(\displaystyle\frac{\tanh M}{\tanh m}\Big)}.
\end{equation}
Hence, we can write
\begin{equation*}
	\tilde{h}(\mu) = \frac{(\beta-\alpha)}{2\ln\Big(\displaystyle\frac{\tanh M}{\tanh m}\Big)}\ln\Big(\displaystyle\frac{\tanh \mu/2}{\tanh m}\Big).
\end{equation*}
We are now in position to compute the Dirichlet energy of the function $\xi_\varepsilon$\,.
\blu{By \eqref{testfootball}, $\big|\nabla_{(x,y,z)} \xi_\varepsilon(x,y,z) \big|=\big|\nabla_{(x,y,z)} h(x,y,z)\big|$, where here we stress the dependence of the nabla operator on the variables. 
By changing variables with \eqref{EllipsHConstr} and by recalling that 
$$(J_{\Psi_a}(\mu,\nu,\varphi))^\top \nabla_{(x,y,z)} h(x,y,z) = \nabla_{(\mu,\nu,\varphi)} \tilde{h}(\mu),$$
owing to \eqref{jacobian}, \eqref{jacobian_determinant}, and \eqref{larichiamiamo}, we obtain that 
$$\nabla_{(x,y,z)} \xi_\varepsilon = \frac{c}{a\sinh\mu(\sin^2\nu + \sinh^2\mu)}
\begin{pmatrix}
- \cosh\mu \sin\nu \cos\varphi\\
 \sinh\mu \cos\nu\\
- \cosh\mu \sin\nu \sin\varphi
\end{pmatrix}.$$
}
Therefore, by using che change of variables theorem and recalling \eqref{jacobian_determinant} and \eqref{eq:BS_elliptic_competitor}, we obtain 
\begin{equation}\label{hUB3}
\begin{split}
&\, \frac{1}{2} \int_{E(a,M) \setminus \overline{E}(a,m)} |\nabla \xi_\eps|^2 \, \de \mathbf{x} \\
=& \, \frac{c^2}{2a^2} \int_0^{2\pi} \!\! \int_0^{\pi} \!\! \int_{2m}^{2M}  \frac{\cosh^2\mu \sin^2\nu +\sinh^2\mu \cos^2\nu}{\sinh^2\mu (\sin^2\nu+\sinh^2\mu)^2} \big|\det J_{\Psi_a}(\mu,\nu,\varphi)\big| \,\de\mu \de \nu \de \varphi \\
=&\, \pi c^2 a \int_0^\pi \!\! \int_{2m}^{2M} \frac{\sin \nu}{\sinh \mu} \,\de\mu \de \nu
=  \frac{\pi a(\beta-\alpha)^2}{\blu{2}\ln \left( \displaystyle\frac{\tanh M}{\tanh m} \right)},
\end{split}
\end{equation}
\blu{where, in the second equality, we have used the elementary relations between trigonometric and hyperbolic functions.}
In the above expression, there are still two choices that we can make: that of the parameters $a$ and $m$.
We now want to choose them in such a way that
\begin{equation}
	\label{fit}
	(N_\eps \cap \{{x = \pm \eps}\})^\circ \subset \overline{E(a,m) \cap \{x = 0 \}}\pm(\eps,0,0).
\end{equation}
To enforce \eqref{fit}, we note that \eqref{EllipsHConstr} implies that, for all $(x,y,z) \in  \partial E(a,m)$ such that $x =0$, it holds
\begin{equation}\label{amconditionellipsoid}
	\frac{y^2}{a^2 \cosh^2 (\blu{2}m)}+\frac{z^2}{a^2 \sinh^2 (\blu{2}m)}   = 1.
\end{equation}
Therefore, choosing $a$ and $m$ to satisfy
\begin{equation}\label{def_a}
a \sinh (\blu{2}m) = \blu{2} \eta,\quad\quad\quad
a \cosh (\blu{2}m) = \blu{2} \delta
\end{equation}
\blu{implies that \eqref{amconditionellipsoid} holds with $\leq$ and therefore}
guarantees the validity of \eqref{fit}.

We now want to get an asymptotic estimate of \eqref{hUB3}, taking into account the fact that all the regimes in this paper consider the case in which $\eta \ll \delta$.
\blu{By \eqref{def_a}, $a^2=4(\delta^2-\eta^2)$, and then
\begin{equation}\label{eq:asympt_a}
a\approx 2\delta\quad \text{ as }\quad \eps\to0.
\end{equation}
} 
Moreover, \blu{from \eqref{def_a},}
\begin{equation}\label{tanh2m}
	\tanh (\blu{2}m) = \frac{\eta}{\delta}.
\end{equation}
Observe that in our regimes, when $\eta \ll \delta$, then $m \ll 1$ and \blu{\eqref{tanh2m} implies that $\tanh(m)\approx\eta/2\delta$ as $\varepsilon\to0$, so that}
\begin{equation}\label{eq23}
	\ln \left( \frac{\tanh M}{\tanh m} \right) = \ln \tanh M - \ln \tanh m \approx  -\ln \tanh m \approx -\ln \frac{\eta}{\blu{2}\delta}\blu{\approx -\ln\frac{\eta}{\delta}},
\end{equation}
for $\eps$ small enough.
This, together with \eqref{hUB3} implies that, for $\eps$ small enough,
\blu{
\begin{equation}\label{footballenergy}
	\lim_{\eps\to0} \frac{\left|\ln (\eta/\delta)\right|}{\delta} \frac{1}{2} \int_{E(a,M) \setminus \overline{E}(a,m)} |\nabla \xi_\eps|^2 \, \de \bold{x} = \pi(\beta-\alpha)^2.
\end{equation}
}

Finally, we note that the elliptic competitor just built gives a better upper bound on the energy of the minimiser $u_\varepsilon$ than the one that could be obtained in \cite{KohSla}, with a spherical harmonic function.
Indeed, in the spherical harmonic case, they obtained an upper estimate with a term of order $\delta$. Therefore, we can notice that
\[
\frac{\delta}{|\ln(\eta/\delta)|}\ll\delta,
\]
as $\eps\to0$. Thus, we obtained a competitor whose order of energy is asymptotically lower than the previous one. This is particularly relevant since such a competitor follows the geometry of our problem, in which the shape of the neck presents the $y$ coordinate way smaller than the $z$ coordinate, ruled by $\delta$ and $\eta$ respectively.

\subsection{Mixed competitor}
The idea now, is to mix the affine competitor in the neck, together with the ellipsoidal just built. The purpose of this new competitor, is to describe whenever the transition happens simultaneously inside and outside the neck.
Consider $A,B\in\R{}$ such that $\alpha\leq A\leq B\leq \beta$.
Let $h\colon E^\ell(a,M) \setminus E^\ell(a,m)\to\R{}$ be the solution to
\begin{equation*}
	\begin{cases}
		\Delta w=0 & \text{in }  E^\ell(a,M) \setminus \overline{E}^\ell(a,m), \\[5pt]
		w = \alpha & \text{on } \partial E^\ell(a,M), \\[5pt]
		w = A & \text{on } \partial E^\ell(a,m),
	\end{cases}
\end{equation*}
and $g\colon E^\rr(a,M) \setminus E^\rr(a,m)\to\R{}$ the solution to
\begin{equation*}
	\begin{cases}
		\Delta w=0 & \text{in } E^\rr(a,M) \setminus \overline{E}^\rr(a,m), \\[5pt]
		w= \beta & \text{on } \partial E^\rr(a,M), \\[5pt]
		w = B & \text{on } \partial E^\rr(a,m);
	\end{cases}
\end{equation*}
\blu{let $h_\eps(x,y,z)\coloneqq h(x+\eps,y,z)$ and $g_\eps(x,y,z)\coloneqq g(x-\eps,y,z)$.}
Recalling \eqref{testfootball} and \eqref{EllipsHConstr}, we define the function $\xi_\varepsilon \colon \R{3}\to\R{}$ as
\begin{equation}\label{mixedcompetitor}
	\xi_\eps(\bfx)\coloneqq
	\begin{cases}
		\alpha & \text{in }\Omega_\eps^\ell\setminus \overline{E}_\eps^\ell(a,M), \\[5pt]
	h_\eps(x,y,z) & \text{in } E_\eps^\ell(a,M)\setminus \overline{E}^\ell(a,m), \\[5pt]
	A & \text{in } E_\eps^\ell(a,\blu{m}), \\[5pt]
		\displaystyle\frac{B-A}{2\eps}x+\frac{B+A}{2} & \text{in } N_\eps, \\[5pt]
	B & \text{in } E_\eps^\rr(a,\blu{m}), \\[5pt]
	g_\eps(x,y,z) & \text{in } E_\eps^\rr(a,M)\setminus \overline{E}^\rr(a,m), \\[5pt]
		\beta & \text{in } \Omega_\eps^\rr\setminus \overline{E}_\eps^\rr(a,M).
	\end{cases}
\end{equation}
\blu{We note that $\xi_\varepsilon$ is an admissible test function for the minimization problem satisfied by~$u_\varepsilon$. Indeed, by the maximum principle, the functions $g$ and $h$ are uniformly bounded in~$\varepsilon$. Therefore, we also get that
\[
\int_{\Omega_\varepsilon} |\xi_\varepsilon - u_{0,\varepsilon}|^2 \de \mathbf{x} \leq C (|E_\varepsilon^l| + |E_\varepsilon^r|),
\]
for $\varepsilon$ small enough.
}
We now want to estimate, asymptotically, its energy.
Using the same argument used to obtain \eqref{hmu}, we can write the explicit solution of the problems above as
\begin{equation*}
		h(\mu) = c^\ell \ln |k^\ell \tanh(\mu/2)|,\quad\quad\quad\text{and}\quad\quad\quad
		g(\mu) = c^\rr \ln |k^\rr \tanh(\mu/2)|.
\end{equation*}
We can explicitly obtain $c^\ell,k^\ell,c^\rr, k^\rr\in\R{}$ by imposing the boundary conditions and we get 
\begin{align*}
	c^\ell=\frac{\alpha-A}{\ln\Big(\displaystyle\frac{\tanh M}{\tanh m}\Big)},\quad\quad\quad
	c^\rr=\frac{\beta-B}{\ln\Big(\displaystyle\frac{\tanh M}{\tanh m}\Big)}
\end{align*}
and 
\begin{align*}
	k^\ell=\frac{\exp\Big(\displaystyle\ln\abs{\frac{\tanh M}{\tanh m}}\frac{\blu{\alpha}}{\alpha-A}\Big)}{\tanh M}\,\quad\quad\quad k^\rr=\frac{\exp\Big(\displaystyle\ln\abs{\frac{\tanh M}{\tanh m}}\frac{\blu{\beta}}{\beta-B}\Big)}{\tanh M}.
\end{align*}
Arguing like in \eqref{hUB3} \blu{and recalling \eqref{eq23}}, we get that 
\begin{align}\label{BoundA}
	\frac{1}{2}\int_{E_\eps^\ell(a,M)\setminus \overline{E}^\ell(a,m)}|\nabla h_\eps|^2\,\de \bold{x}=\frac{\pi a(A-\alpha)^2}{\abs{\ln(\eta/\delta)}}
\end{align}
and 
\begin{align}\label{BoundB}
\frac{1}{2}\int_{E_\eps^\rr(a,M)\setminus \overline{E}^\rr(a,m)}|\nabla g_\eps|^2\ \de \bold{x}=	\frac{\pi a(B-\beta)^2}{\abs{\ln(\eta/\delta)}},
\end{align}
as $\eps\to0$.
Therefore, from \eqref{AffineBound}, \eqref{BoundA}, and \eqref{BoundB}, we obtain
\begin{align*}
	F(\xi_\eps,\Omega_\eps)= \frac{\pi a}{|\ln(\eta/\delta)|}\big[(A-\alpha)^2+(B-\beta)^2\big]+\frac{\delta \eta}{\eps}(B-A)^2.
\end{align*}
Recalling \eqref{eq:asympt_a}, for $\eps$ small enough, we can write
\begin{align}\label{lowerbound}
	F(\xi_\eps,\Omega_\eps)\approx \frac{\blu{2}\pi \delta}{|\ln(\eta/\delta)|}\big[(A-\alpha)^2+(B-\beta)^2\big]+\frac{\delta \eta}{\eps}(B-A)^2.
\end{align}
Now we compute the minimum of the right-hand of \eqref{lowerbound} under the constraint that $\alpha\leq A\leq B\leq\beta$.  It is possible to see that a solution is in the interior of the admissible region, and thus the optimal $A$ and $B$ are given by the solution of the system
\begin{align*}
	\begin{cases}
		\displaystyle\frac{\blu{2}\pi}{|\ln(\eta/\delta)|}(A-\alpha)-\frac{ \eta}{\eps}(B-A)=0,\\[10pt]
\displaystyle	\frac{\blu{2}\pi }{|\ln(\eta/\delta)|}(B-\beta)+\frac{ \eta}{\eps}(B-A)=0,
	\end{cases}
\end{align*} 
which are
\begin{align}\label{optimalAandB}
	 A=\frac{\displaystyle\frac{\pi\alpha}{\abs{\ln(\eta/\delta)}}+\frac{\eta}{\eps}\blu{\frac{\alpha+\beta}{2}}}{\displaystyle\frac{\pi}{\abs{\ln(\eta/\delta)}}+\frac{\eta}{\eps}}\quad\text{and}\quad
	B=\frac{\displaystyle\frac{\pi\beta}{\abs{\ln(\eta/\delta)}}+\frac{\eta}{\eps}\blu{\frac{\alpha+\beta}{2}}}{\displaystyle\frac{\pi}{\abs{\ln(\eta/\delta)}}+\frac{\eta}{\eps}}.
\end{align}
The explicit values of $A$ and $B$ in \eqref{optimalAandB} will be crucial in the various regimes when we will need to infer the boundary conditions of the rescaled profile at the edge of the neck.
In conclusion, from \eqref{lowerbound}, if $u_\eps$ is a local minimiser, we then have
\begin{align}\label{upperboundueps}
	F(u_\eps,\Omega_\eps) & \leq\frac{\blu{2}\pi\delta}{|\ln(\eta/\delta)|}\big[(A-\alpha)^2+(B-\beta)^2\big]+\frac{\delta \eta}{\eps}(B-A)^2.
\end{align}
Finally, notice that the right-hand side of \eqref{upperboundueps} has a clear separation between the energetic contribution of the competitor inside and outside the neck, as well as their orders of the energy.

\section{Analysis of the problem in the several regimes}
\label{sec:analysis}
\indent

In this section we carry out the rigorous analysis of the asymptotic behaviour of the solution, obtaining information on its energy and its behaviour inside and close to the neck.
To this aim, define 
\[
N\coloneqq [-1,1]^3,
\]
which is the neck $N_\eps$ rescaled to size of order $1$, that is, under the change of coordinates $(x,y,z)\to(x/\eps,y/\delta,z/\eta)$.
\blu{In the following subsections, we will perform various rescalings and we will always denote by $\widetilde{\Omega}_\varepsilon$ the corresponding rescaled domains.}

\subsection{Super-thin neck}

In this regime the parameters are ordered as
 $\eps \gg \delta \gg \eta$. Namely, we have
 \begin{align*}
 	\lim_{\eps\to0} \frac{\delta}{\eps}=0\quad\text{and}\quad\lim_{\eps\to0}\frac{\eta}{\delta}=0.
 \end{align*}
According to the heuristics in Section \ref{sec:heuristics}, we expect the transition to happen entirely inside the neck.
If $u_\eps$ is a local minimiser of the functional \eqref{004}, the \blu{convenient rescaling} is
\begin{align*}
	v_\eps(x,y,z)\coloneqq u_\eps(\eps x,\delta y, \eta z).
\end{align*}
Define $\widetilde{\O}_\eps$, $\widetilde{\O}^\ell_\eps$, and $\widetilde{\O}^\rr_\eps$ as the rescaled domain $\O_\eps$, $\O^\ell_\eps$, and $\O^\rr_\eps$, respectively.
Note that, as $\eps\to0$, \blu{the closure of the set $\widetilde{\O}_\eps$ converges locally in the Hausdorff sense (see Definition \ref{def:Hausdorff}) to the closure of the set}
\[
\Omega_\infty \coloneqq \Omega^\ell_\infty\cup N \cup \Omega^\rr_\infty,
\]
where $\Omega^\ell_\infty\coloneqq \{x<-1\}$ and $\Omega^\rr_\infty \coloneqq \{x>1\}$.
\blu{The strategy of the proof of Theorem \ref{Theoremsuperthin} and of Theorem \ref{Theoremflatthin} is similar to that employed in \cite[Theorem 4.1]{KohSla}. For clarity, we chose to present it in full details also in here.}

\begin{theorem}\label{Theoremsuperthin}
Let $(u_\eps)_\eps$ be an admissible family of local minimisers as in Definition \ref{def:loc_min}.
Assume $\eps \gg \delta \gg \eta$.
Then,
	\begin{align*}
		\lim_{\eps\to0} \frac{\eps}{\delta\eta}F(u_\eps,\O_\eps)
		=\lim_{\eps\to0} \frac{\eps}{\delta\eta}F(u_\eps, N_\eps)
		= (\beta-\alpha)^2.
	\end{align*}
Moreover, for $\varepsilon>0$ let $v_\varepsilon \colon \widetilde{\O}_\eps\to\R{}$ be defined as
\[
v_\eps(x,y,z)\coloneqq u_\eps(\eps x,\delta y, \eta z).
\]
Then, the following hold:
\begin{itemize}
\item[(i)] $v_\eps\chi_{\widetilde{\O}^\ell_\eps}- \alpha\chi_{\widetilde{\O}^\ell_\eps}\to0$ in $L^6(\R{3})$ and 
$v_\eps\chi_{\widetilde{\O}^\rr_\eps}- \beta\chi_{\widetilde{\O}^\rr_\eps}\to0$ in $L^6(\R{3})$, as $\eps\to0$;
\item[(ii)] There exists $\hat{v}\in H^1(N)$ such that $v_\eps\wto \hat{v}$ weakly in $H^1(N)$, as $\eps\to0$;
\item[(iii)] It holds that $\hat{v}(x,y,z)=v(x)$, where $v\in H^1(-1,1)$ is the unique minimiser of the variational problem
	\begin{align*}
		\min\bigg\{ \frac{1}{2} \int_{-1}^1 |v'|^2\ \de x:\ v\in H^1(-1,1),\ v(-1)=\alpha,\ v(1)=\beta\bigg\}.
	\end{align*}
In particular, $v(x)=\dfrac{\beta-\alpha}{2}x+\dfrac{\alpha+\beta}{2}$.
\end{itemize}
\end{theorem}

\begin{proof}
\emph{Step 1: existence of $\hat{v}$.} 
Note that, using assumption (H3), for $\eps>0$ sufficiently small, it holds that
\[
2N\cap \O_\infty\subset \widetilde{\O}_\eps.
\]
The reason why we take $2N$ and not only $N$ is because we need the boundary conditions to converge, \blu{and this is an easy way to ensure such a convergence}.
We claim that
\[
\sup_{\varepsilon>0} \|v_\varepsilon\|_{H^1(2N\cap \O_\infty)} < \infty.
\]
First of all, using the fact that $\eps \gg \delta \gg \eta$, we get that
\begin{equation}\label{eq:est_grad_v_eps}
\begin{split}
\|\nabla v_\varepsilon\|^2_{L^2(2N\cap \O_\infty)}=
&\frac{\varepsilon}{\delta\eta} \int_{2N_\varepsilon\cap \O_\varepsilon} \Big((\partial_x u_\eps)^2 + \frac{\delta^2}{\eps^2}(\partial_y u_\eps)^2 + \frac{\eta^2}{\eps^2}(\partial_z u_\eps)^2\Big)\,\de\textbf{x} \\
\leq&\frac{\varepsilon}{\delta\eta} F(u_\varepsilon,\O_\eps)
\leq C <\infty,
\end{split}
\end{equation}
where the first step follows by using a change of variable, while the second from \eqref{AffineBound} with $A=\alpha$ and $B=\beta$.
Moreover, 
\[
\int_{2N\cap \O_\infty} |v_\eps|^2\,\de \mathbf{x}
=\frac{1}{\eps\delta\eta}\int_{2N_\eps\cap\Omega_\eps} |u_\eps|^2\,\de\mathbf{x}
\leq \frac{\big(\sup_{\mathbf{x}\in 2N_\eps\cap \Omega_\eps}|u_\eps|^2\big) |2N_\eps\cap \Omega_\eps|}{\eps\delta\eta}=C.
\]
Thus, we get that, up to a subsequence, $v_\varepsilon$ converges weakly in $H^1(N)$ to a function $\hat{v}\in H^1(N)$.
The independence of the subsequence will follow from Step $2$, where we show that the limit is the \emph{unique} solution to a variational problem.

\noindent\emph{Step 2: limiting problem and behavior inside the neck.}
We now want to characterize the function $v$ as the unique solution to a variational problem.
We do this in two steps: first, we identify a functional that will be minimized, and then we identify the boundary conditions.\\
We have that 
\begin{align}\label{eq:est_energy_stn}
	\liminf_{\eps\to0}\frac{\eps}{\delta\eta} F(u_\eps,\O_\eps)&\geq	\liminf_{\eps\to0}\frac{\eps}{\delta\eta}\Big(\frac{1}{2}\int_{N_\eps}|\nabla u_\eps|^2\,\de \textbf{x} +\int_{N_\eps} W(u_\eps)\, \de\bfx \Big)\nonumber \\[5pt]
	&=\liminf_{\eps\to0}\frac{1}{2}\int_{N} \Big((\partial_x v_\eps)^2 + \frac{\eps^2}{\delta^2}(\partial_y v_\eps)^2 + \frac{\eps^2}{\eta^2}(\partial_z v_\eps)^2\Big)\,\de \textbf{x} +\eps^2\int_N W(v_\eps)\, \de\bfx \nonumber  \\[5pt]
	&\geq\liminf_{\eps\to0}\frac{1}{2}\int_N |\nabla v_\eps|^2\, \de \textbf{x}  
	\geq \frac{1}{2}\int_N |\nabla v|^2\, \de \textbf{x},
	\end{align}
where in the last step we used the fact that $\eps\gg\delta\gg\eta$.
Notice that, from the bound \eqref{eq:est_grad_v_eps} and the fact that $\eps/\eta\to\infty$ and $\eps/\delta\to\infty$, we necessarily have that $v$ does not depend on $y$ and $z$. Namely, $v_\eps$ converges to a function $\hat{v}\in H^1(N)$, of the form $\hat{v}(x,y,z)=v(x)$, where $v\in H^1(-1,1)$. Therefore, from \eqref{eq:est_energy_stn}, we can write
\begin{align}\label{eq:est_energy_stnpt2}
\frac{1}{2}\int_N |\nabla v|^2\, \de \bold{x}=  2\int_{-1}^1 (\hat{v}')^2\ \de x 
\geq&\, 2\min\Big\{ \int_{-1}^1 (w')^2\, \de x:\ w\in H^1(-1,1),\ w(\pm1)=\hat{v}(\pm1)\Big\} \nonumber  \\[5pt]
=&\,(\hat{v}(1)-\hat{v}(-1)\big)^2,
\end{align}
where last step follows by an explicit minimisation.

Now we claim that $\hat{v}(-1)=\alpha$ and $\hat{v}(1)=\beta$.
We prove the former, since the latter follows by using a similar argument.
The idea (introduced in \cite{KohSla}) is to use the scale-invariant Poincar\'{e} inequality
\begin{align}\label{poincare}
\Big(\int_{\Omega_\eps^\ell}|u_\eps-\bar{u}_\eps|^6\, \de \bold{x}\Big)^{\frac{1}{6}}\leq C \Big(\int_{\Omega_\eps^\ell}|\nabla u_\eps|^2\,\de \bold{x}\Big)^\frac{1}{2},
\end{align}
where $C>0$ and $\bar{u}_\eps$ is the average on $\Omega_\eps^\ell$ of $u_\eps$.
Using a change of variable, we estimate (by neglecting the potential term as in \eqref{eq:est_energy_stn}) the right-hand side of \eqref{poincare} as
\begin{align*}
	\eps\delta\eta\int_{\widetilde{\O}_\eps^\ell}|v_\eps-\bar{v}_\eps|^6\,\de \bold{x}
&\leq C \Big(\frac{\delta\eta}{\eps}
	\int_{\widetilde{\O}_\eps^\ell} (\partial_x v_\eps)^2
	+\frac{\eps^2}{\delta^2}(\partial_y v_\eps)^2
	+\frac{\eps^2}{\eta^2}(\partial_z v_\eps)^2\,\de \bold{x}
	\Big)^3\\
&\leq C \Big(\frac{\delta\eta}{\eps}\Bigr)^3
		\left(\frac{\eps}{\delta\eta}F(u_\eps,\O_\eps)\right)^3.
\end{align*}
Now, using the fact that $\delta\eta/\eps^2\to0$ as $\eps\to0$, \blu{using} \eqref{eq:est_grad_v_eps}, we get that 
\begin{align}\label{poincareestimate}
	\int_{\widetilde{\O}_\eps^\ell}|v_\eps-\bar{v}_\eps|^6\,\de \bold{x}\leq C \Big(\frac{\delta\eta}{\eps^2}\Big)^2\left(\frac{\eps}{\delta\eta}F(u_\eps,\O_\eps)\right)^3\to0,
\end{align}
as $\varepsilon\to0$. 
Therefore, $v_\varepsilon\chi_{\widetilde{\O}_\eps^\ell} - \bar{v}_\eps\chi_{\widetilde{\O}_\eps^\ell}\to 0$ in $L^6(\R{3})$, as $\eps\to0$.
We now claim that $\bar{v}_\eps\chi_{\widetilde{\O}_\eps^\ell} - \alpha\chi_{\widetilde{\O}_\eps^\ell}\to 0$ in $L^6(\R{3})$, as $\eps\to0$. Indeed, by definition of local minimiser we have that
\[
||u_\eps-\alpha||_{L^1(\Omega_\eps^\ell)}\to0,
\]
as $\eps\to0$.
Therefore, $\bar{u}_\eps\chi_{\O_\eps^\ell} - \alpha\chi_{\O_\eps^\ell}\to 0$ in $L^6(\R{3})$, which yields that $\bar{v}_\eps\chi_{\widetilde{\O}_\eps^\ell} - \alpha\chi_{\widetilde{\O}_\eps^\ell}\to 0$ in $L^6(\R{3})$.
Thus, since $v_\eps\to\hat{v}$ strongly in $L^2(2N\cap\O_\infty)$ as $\eps\to0$, we get that $\hat{v}(-1)=\alpha$.

\noindent\emph{Step 3: asymptotic behaviour of the energy.}
From \eqref{eq:est_energy_stn} and \eqref{eq:est_energy_stnpt2} we can conclude that
\begin{align}\label{liminfsuperthin}
	\liminf_{\eps\to0}\frac{\eps}{\delta\eta} F(u_\eps,\O_\eps)\geq (\beta-\alpha)^2.
\end{align}
On the other hand, denoting by $\xi_\eps$ the affine competitor in \eqref{affinecompetitor}, we have that
\begin{align}\label{limsupsupethin}
\limsup_{\eps\to 0}\frac{\eps}{\delta\eta} F(u_\eps,\O_\eps)\leq\limsup_{\eps\to 0} \frac{\eps}{\delta\eta}F(\xi_\eps,N_\eps)=(\beta-\alpha)^2.
\end{align}
Thus, from  \eqref{liminfsuperthin}, and \eqref{limsupsupethin}, we get
\begin{align*}
	\lim_{\eps\to0} \frac{\eps}{\delta \eta} F(u_\eps,\O_\eps)=(\beta-\alpha)^2.
\end{align*}
In particular, we get that all inequalities in \eqref{eq:est_energy_stn} are equalities, proving that
\[
\lim_{\eps\to0} \frac{\eps}{\delta\eta}F(u_\eps,\O_\eps)
		=\lim_{\eps\to0} \frac{\eps}{\delta\eta}F(u_\eps, N_\eps).
\]
This concludes the proof.
\end{proof}

\subsection{Flat-thin neck}

In this regime the parameters are ordered as
$\eps \approx\delta \gg \eta$. Namely, we have
\begin{align*}
	\lim_{\eps\to0} \frac{\delta}{\eps}=m\in(0,+\infty)\quad\text{and}\quad \lim_{\eps\to 0}\frac{\eta}{\eps}=0.
\end{align*}
In this case, the behaviour of an admissible family of local minimisers is similar to the super-thin regime and strategy of the proof is similar to that of Theorem \ref{Theoremsuperthin}. Therefore, we only highlight the main differences. Since we expect the transition to happen entirely inside the neck, we would like to use a rescaling for which the neck $N_\eps$ transforms in $N\coloneqq [-1,1]^3$ \blu{(where we take, without loss of generality, $m=1$)}. Given a local minimizer $u_\eps$ of the functional~\eqref{004}, the \blu{convenient rescaling} is
\begin{align*}
	v_\eps(x,y,z)\coloneqq u_\eps(\eps x,\eps y, \eta z).
\end{align*}
If we rescale in this way, the limiting domain \blu{(in the sense of local Hausdorff convergence)} becomes
$$\Omega_\infty = \Omega^\ell_\infty \cup N \cup \Omega^\rr_\infty,$$
where $\Omega^\ell_\infty = \{x<-1\}$ and $\Omega^\rr_\infty = \{x>1\}$. 

\begin{theorem}\label{Theoremflatthin}
Let $(u_\eps)_\eps$ be an admissible family of local minimisers as in Definition \ref{def:loc_min}.
	Assume $\eps \approx \delta \gg \eta$. Then,
	\begin{align*}
		\lim_{\eps\to0} \frac{1}{\eta}F(u_\eps,\O_\eps)
		=\lim_{\eps\to 0} \frac{1}{ \eta} F(u_\eps,N_\eps)
		= (\beta-\alpha)^2.
	\end{align*}
Moreover, for $\varepsilon>0$ let $v_\varepsilon \colon \widetilde{\O}_\eps\to\R{}$ be defined as
\[
v_\eps(x,y,z)\coloneqq u_\eps(\eps x,\eps y, \eta z).
\]
Then, the following hold:
\begin{itemize}
\item[(i)] $v_\eps\chi_{\widetilde{\O}^\ell_\eps}- \alpha\chi_{\widetilde{\O}^\ell_\eps}\to0$ in $L^6(\R{3})$ and 
$v_\eps\chi_{\widetilde{\O}^\rr_\eps}- \beta\chi_{\widetilde{\O}^\rr_\eps}\to0$ in $L^6(\R{3})$ as $\varepsilon\to0$;
\item[(ii)] There exists a function $\hat{v}\in H^1(N)$ such that $v_\eps\wto \hat{v}$ weakly in $H^1(N)$ as $\varepsilon\to0$;
\item[(iii)] It holds that $\hat{v}(x,y,z)=v(x)$, where $v\in H^1(-1,1)$ is the unique minimiser of the variational problem
	\begin{align*}
		\min\bigg\{ \frac{1}{2} \int_{-1}^1 |v'|^2\,\de x:\ v\in H^1(-1,1),\ v(-1)=\alpha,\ v(1)=\beta\bigg\}.
	\end{align*}
In particular, $v(x)=\dfrac{\beta-\alpha}{2}x+\dfrac{\alpha+\beta}{2}$.
\end{itemize}
\end{theorem}

\begin{proof}
In the same way as in Theorem \ref{Theoremsuperthin}, we can prove that
\begin{align}\label{boundFflatthin}
	\sup_{\eps>0}\|\nabla v_\eps\|^2_{L^2(N)}\leq \frac{\eps}{\delta\eta}F (u_\eps,\O_\eps)\leq C<\infty.
\end{align}
and
\[
\sup_{\eps>0} \|v_\eps\|_{L^2(N)}<C.
\]
Therefore, by compactness there exists $v\in H^1(N)$ such that, up to a subsequence,  $v_\eps\wto v$ in $H^1(N)$. The independence of the subsequence will follow from the fact that the limit is the \emph{unique} solution to a variational problem.  \\
	Thus,  we can write
	\begin{align}\label{liminfflatthin}
	\nonumber	\liminf_{\eps\to0}\frac{1}{\eta} F(u_\eps,N_\eps)&= \liminf_{\eps\to0}\frac{1}{2} \int_{N} \Big((\partial_x v_\eps)^2 + (\partial_y v_\eps)^2 + \frac{\eps^2}{\eta^2}(\partial_z v_\eps)^2\Big)\,\de\mathbf{x}+\eps^2\eta\int_{N} W(v_\eps)\,\de\mathbf{x}\\[5pt]
	\nonumber&\geq\int_{[-1,1]^2} \big((\partial_x\hat{v})^2+(\partial_y \hat{v})^2\big)\, \de x\de y\\[5pt]
	\nonumber&\geq\min\bigg\{ \int_{[-1,1]^2} \big((\partial_x w)^2+(\partial_y w)^2\big)\, \de x\de y:\ w\in H^1([-1,1]^2),\\[5pt] 
		& \phantom{\geq\min \bigg\{\;} w(\pm1,y)=\hat{v}(\pm1,y),\ \forall y\in[-1,1]\bigg\},
	\end{align}
	where in the previous to last step we used \eqref{boundFflatthin} and the fact that $\eps/\eta\to+\infty$. Then we necessarily have that $v$ does not depend on $z$. Namely, $v_\eps$ converges to a function $\hat{v}\in H^1(N)$, of the form $\hat{v}(x,y,z)=v(x,y)$, where $v\in H^1([-1,1]^2)$. \\
	Now, would like to show that the boundary conditions $\hat{v}(\pm1,y)$ are independent from $y$. This is done by acting similarly as in \eqref{poincare} and \eqref{poincareestimate}. Indeed by using the scale-invariant Poincar\'{e} inequality \eqref{poincare}, we get
\begin{align*}
\int_{\Omega_\eps^\ell}|v_\eps-\bar{v}_\eps|^6\, \de \bold{x}\leq C \Big(\frac{\eta}{\eps}\Big)^2\Big(\frac{1}{\eta}F(v_\eps,\O_\eps)\Big)^3.
\end{align*}
From \eqref{boundFflatthin} and the fact that $\eta/\eps\to0$ as $\eps\to0$ we can conclude that $v=\alpha$ on $\Omega^\ell_\infty$ and $v=\beta$ on $\Omega^r_\infty$\,, independently on $y$. Therefore, from \eqref{liminfflatthin} we get
\begin{align}\label{eq:est_Neps}
	\liminf_{\eps\to0}\frac{1}{\eta} F(u_\eps,N_\eps)
	&\geq\min\bigg\{ \int_{[-1,1]^2} \big((\partial_x w)^2+(\partial_y w)^2\big)\, \de x\de y:\ w\in H^1([-1,1]^2),\nonumber \\[5pt] 
	&\phantom{\geq\min\bigg\{} w(-1,y)=\alpha,\ w(1,y)=\beta,\ \forall y\in[-1,1]\bigg\} \nonumber \\[5pt]
	&\geq\min\bigg\{ \int_{[-1,1]^2} (\partial_x w)^2\ \de x\de y:\ w\in H^1([-1,1]^2),\nonumber \\[5pt] 
	&\phantom{\geq\min\bigg\{} w(-1,y)=\alpha,\ w(1,y)=\beta,\ \forall y\in[-1,1]\bigg\}\nonumber \\[5pt]
	&=(\beta-\alpha)^2,
\end{align}
where the last step follows by an explicit computation.
On the other hand, denoting by $\xi_\eps$ the affine competitor in \eqref{affinecompetitor}, we have that
\begin{align}\label{limsupsupethin_1}
	\limsup_{\eps\to 0}\frac{1}{\eta} F(u_\eps,\O_\eps)\leq\limsup_{\eps\to 0} \frac{1}{\eta}F(\xi_\eps,N_\eps)=(\beta-\alpha)^2.
\end{align}
Finally, using \eqref{eq:est_Neps} and \eqref{limsupsupethin_1}, we get
\begin{align*}
	\lim_{\eps\to0} \frac{1}{ \eta} F(u_\eps,\O_\eps)=\lim_{\eps\to 0} \frac{1}{ \eta} F(u_\eps,N_\eps)=(\beta-\alpha)^2.
\end{align*}
This concludes the proof.
\end{proof}

\subsection{Interlude: convergence of Neumann problems}

In this short interlude, we recall a convergence result for solutions to elliptic problems with Neumann boundary conditions that will be crucial to carry out the analysis of the asymptotic behaviour of the rescaled profiles outside the neck.
\blu{We first recall the notion of Hausdorff convergence.}

\begin{defin}\label{def:Hausdorff}
We say that a sequence of closed sets $(A_n)_n\subset\R{3}$ \emph{converges in the Hausdoff metric} to a closed set $A$ if
\[
\lim_{n\to\infty} \max\left\{ \sup_{\bold{x}\in A}\mathrm{d}(\bold{x}, A_n),\,
	\sup_{\bold{y}\in A_n}\mathrm{d}(\bold{y}, A)  \right\} =0.
\]
Here, $\mathrm{d}(\bold{x}, A)$ denotes the distance between the point $\bold{x}\in\R{3}$ and the set $A$.
We denote this convergence by $A_n\stackrel{\mathrm{H}}{\rightarrow}A$.
We say that the convergence is local if it holds on every compact set, namely if $A_n\cap K \stackrel{\mathrm{H}}{\rightarrow} A\cap K$, for every compact set $K$.
\end{defin}

The result is the following. For a proof, we refer to \cite[Proposition 6.2]{MorSla12} (see also \cite{ChaDov, DMEboPon}).
Note that the argument in there is only detailed for the case of $\R{2}$. Nevertheless, all of the computations carry out also in the three-dimensional case, including the $W^{2,p}_{\loc}$ convergence, which follows from standard interior regularity estimates.

\begin{theorem}\label{neumanconvergence}
	Let $(\O_\eps)_\eps\subset\R{3}$ be a sequence of bounded open sets \blu{with Lipschitz boundary} such that, as $\eps\to0$,
	\begin{align*}
	\chi_{\O_\eps}\to\chi_{\O_\infty}\,\, \text{ in } L^1_{\loc}(\R{3})\qquad \text{ and }\qquad
	\R{3}\setminus \O_\eps \stackrel{\mathrm{H}}{\to} \R{3}\setminus\O_\infty
	\quad\text{ locally},
\end{align*}
for some open set $\O_\infty$\,. Moreover, assume that $\O_\eps\subset\O_\infty$\,, for every $\eps>0$.
Let $p > 2$, and let $(f_\eps)_\eps\subset L^p_{\loc}(\R{3})$ be such that, as $\eps\to0$,
\[
f_\eps\chi_{\O_\eps}\to f_\infty\chi_{\O_\infty} \quad\text{ in } L^p_{\loc}(\R{3}),
\]
for some $f_\infty\in L^p_{\loc}(\R{3})$.
Let $u_\eps\in H^1(\O_\eps)$ be a weak solution to
\[
\left\{
\begin{array}{ll}
\Delta u_\eps = f_\eps &\text{ in } \O_\eps,\\
\partial_\nu u_\eps = 0 &\text{ on } \partial\O_\eps.
\end{array}
\right.
\]
Assume that $(u_\eps\chi_{\O_\eps})_\eps$ is locally equi-bounded in $L^\infty(\R{3})$.
Then, up to a subsequence, 
\[
u_\eps\chi_{\O_\eps}
	\to \hat{u}\chi_{\O_\infty}\quad 
	\text{in } L^q_{\loc}(\mathbb{R}^3) \text{ for all } q\in[1,\infty),
\]
as $\eps\to0$, and
\[
\nabla u_\eps\chi_{\O_\eps}
	\to \nabla \hat{u}\chi_{\O_\infty}\quad
	\text{in } L^2_{\loc}(\mathbb{R}^3;\mathbb{R}^3),
\]
as $\eps\to0$, where $\hat{u}\in H^{1}_{\loc}(\O_\infty)$ is a weak solution to
\[
\left\{
\begin{array}{ll}
\Delta \hat{u} = f_\infty &\text{ in } \O_\infty,\\
\partial_\nu \hat{u} = 0 &\text{ on } \partial\O_\infty.
\end{array}
\right.
\]
Moreover, $u_\eps\to \hat{u}$ in $W^{2,p}_{\loc}(\O_\infty)$, as $\eps\to0$.
\end{theorem}

\begin{remark}
The reason why in the above result the convergence of the complements of the open sets $\O_\eps$ is required, and the fact that the limiting set has to be open, is in order to ensure that at each point of the limiting set there is only one side where the limiting set is.
For instance, we want to avoid situations of the type
\[
\O_\eps\coloneqq \{ (\cos\theta, \sin\theta) : \theta\in (0,2\pi-\eps) \},
\]
or of the type
\[
\O_\eps\coloneqq (0,1)^2 \cup \left( [1,2)\times (-\eps,\eps) \right).
\]
In both cases, the limiting set has part of the topological boundary that creates troubles in defining the limiting PDE.
\end{remark}

\subsection{Window thick regime}

Here we consider the scaling of the energy in the window thick regime $\delta\gg\eta\gg\eps$, namely where
\begin{align*}
\lim_{\eps\to0}\frac{\eta}{\delta}=0\quad\text{and}\quad\lim_{\eps\to0} \frac{\eta}{\eps}=+\infty.
\end{align*}
Since the ellipsoidal competitor outside the neck provides an energy whose order is lower than the energy of the affine competitor in the neck, we expect the transition happening outside the neck. 
If $u_\eps$ is a local minimiser of the functional \eqref{004}, 
the \blu{convenient rescaling} that allows us to both see a nice limiting space, and to use the scale-invariant Poincar\'{e} inequality \eqref{poincare} is
\begin{align*}
	v_\eps(x,y,z)\coloneqq u_\eps(\delta x, \delta y, \delta z).
\end{align*}
Using this rescaling, the limiting domain \blu{(in the sense of local Hausdorff convergence)} becomes
\[
\widetilde{\Omega}_\infty = \{ x<0 \} \cup ( \{0\}\times [-1,1] \times \{0\} ) \cup \{ x>0 \}.
\]
However, note that
\[
\R{3}\setminus \widetilde{\Omega}_\varepsilon
\stackrel{\mathrm{H}}{\rightarrow} 
\R{3}\setminus \Omega_\infty\,,
\]
where
\[
\widetilde{\Omega}_\varepsilon \coloneqq \frac{1}{\delta}\Omega_\varepsilon\,,\qquad
\Omega_\infty \coloneqq \{ x<0 \} \cup \{ x>0 \}.
\]
We are now in position to prove the main result of this section.

\begin{theorem}\label{Theoremwindowthick}
Let $(u_\eps)_\eps$ be an admissible family of local minimisers as in Definition \ref{def:loc_min}.
Assume that $\delta\gg\eta\gg\eps$.
Then,
	\begin{equation*}
		\lim_{\eps\to0} \frac{|\ln(\eta/\delta)|}{\delta}F(u_\eps,\O_\eps)=
		\lim_{\eps\to0} \frac{|\ln(\eta/\delta)|}{\delta}F(u_\eps, \O_\eps\setminus N_\eps)=
		 \pi(\beta-\alpha)^2. 
	\end{equation*}
Moreover, for $\varepsilon>0$ let $v_\varepsilon \colon \widetilde{\O}_\eps\to\R{}$ be defined as
\[
v_\eps(x,y,z)\coloneqq u_\eps(\delta x, \delta y, \delta z).
\]
Then, the following statements hold:
\begin{itemize}
\item[(i)] $v_\eps\to \frac{\alpha+\beta}{2}$ uniformly on $B_R$, as $\eps\to0$, \blu{for all $2R<r_0$, where $r_0$ is given by (H3)};
\item[(ii)] There exists $\hat{v}\in H^1_{\loc}(\Omega_\infty)$ such that $v_\eps\chi_{\widetilde{\O}_\eps}\to \hat{v}\chi_{\O_\infty}$ strongly in $H^1_{\loc}(\R{3})$, as $\eps\to0$;
\item[(iii)] The function $\hat{v}$ is the unique minimiser of the variational problem
	\begin{align*}
		\min\bigg\{ \frac{1}{2} \int_{\O_\infty} |\nabla v|^2\, \de x:
		v\in \mathcal{A}  \bigg\},
	\end{align*}
where
\[
\mathcal{A}\coloneqq \bigl\{v\in H^1_{\loc}(\O_\infty),\ v-\alpha\chi_{\O_\infty^\ell}-\beta\chi_{\O^\rr_\infty} \in L^6(\O_\infty),\, v=\frac{\alpha+\beta}{2} \text{ on } B_R\cap\O_\infty \Bigr\},
\]
and $\O^\ell_\infty\coloneqq \{x<0\}$, and $\O^\rr_\infty\coloneqq \{x>0\}$.
\end{itemize}
\end{theorem}

\begin{proof}
\emph{Step 1: lower bound of the energy outside the neck.}
We claim that
\begin{equation}\label{eq:window_lower_bound_energy}
\liminf_{\eps\to0} \frac{|\ln(\eta/\delta)|}{\delta}F(u_\eps,\O_\eps\setminus N_\varepsilon)\geq \pi(\beta-\alpha)^2. 
\end{equation}
The strategy to prove it is the following.
We fix $\rho<R$, and take vanishing sequences $(a_\varepsilon)_\varepsilon$ and $(M_\varepsilon)_\varepsilon$ such that
\begin{equation}\label{eq:a_e-M_e}
a_\eps \sinh(2 M_\eps) = 2 \eta,\quad\quad\quad
a_\eps \cosh(2 M_\eps) = 2 \delta.
\end{equation}
Note that this choice ensures that
\begin{align*}
(N_\eps \cap\{x=-\varepsilon\})^\circ \subset \overline{E^\ell_\eps(a_\eps, M_\eps)\cap\{x=-\varepsilon\}}, \quad
(N_\eps \cap\{x=\varepsilon\})^\circ \subset \overline{E^\rr_\eps(a_\eps, M_\eps)\cap\{x=\varepsilon\}},
\end{align*}
where we define
\begin{align*}
	E^\ell_\eps(a_\eps,M_\eps)&\coloneqq
	\left( E(a_\eps,M_\eps)\cap\{x<0\} \right) + (0,0,-\eps),\\[5pt]
	E^\rr_\eps(a_\eps,M_\eps)&\coloneqq \left(E(a_\eps,M_\eps)\cap\{x>0\}\right) + (0,0,\eps).
\end{align*}

In Step 1.1, we prove that, for all $\theta>0$, there exists $\varepsilon_0>0$ such that
\begin{align}\label{firstboundarycondition_1}
	m_1-\theta &\leq u_\eps \leq m_1+\theta
	\quad\quad\text{on } \overline{E}^\ell_\eps(a_\eps, M_\eps),\\[5pt]
\label{firstboundarycondition_2}
	m_2-\theta &\leq u_\eps \leq m_2+\theta
	\quad\quad\text{on } \overline{E}^\rr_\eps(a_\eps, M_\eps),
\end{align}
for all $\varepsilon<\varepsilon_0$.
In Step 1.2, we prove that, for all $\gamma>0$, there exists $\varepsilon_0>0$ such that
\begin{align}
	\label{secondboundarycondition_1}
	\alpha-\gamma &\leq u_\epsilon \leq \alpha+\gamma
	\quad\quad\text{on }\partial E^\ell(a_\eps, \rho),\\[5pt] 
\label{secondboundarycondition_2}
	\beta-\gamma &\leq u_\eps \leq \beta+\gamma
	\quad\quad\text{on }\partial E^\rr(a_\eps, \rho),
\end{align}
for all $\varepsilon<\varepsilon_0$.
In Step 1.3, we then use this information to first prove that there exists $m_1, m_2\in\R{}$ such that
\begin{equation}\label{eq:window_lower_bound_energy}
\liminf_{\eps\to0} \frac{|\ln(\eta/\delta)|}{\delta}F(u_\eps,\O_\eps\setminus N_\varepsilon)\geq 2\pi\big[  (m_1-\alpha)^2+ (m_2-\beta)^2\big].
\end{equation}
Finally, we obtain the claim of the step by optimizing in $m_1$ and $m_2$ on the right-hand side of \eqref{eq:window_lower_bound_energy}.

\noindent\emph{Step 1.1: boundary conditions on the internal ellipsoid.}
Here, we want to prove the validity of \eqref{firstboundarycondition_1} and of \eqref{firstboundarycondition_2}.
We want to understand the limiting behaviour of $v_\eps$.
The idea is to obtain such information by looking at the limit of the PDE satisfied by the limit of the sequence $(v_\eps)_\eps$.
First, we notice that since $u_\eps$ is a critical point of the energy
$F_\eps$, we have that $u_\eps$ satisfies the Euler-Lagrange equation
\[
	\int_{\Omega_\eps}\nabla u_\eps \cdot \nabla\varphi\, \de\textbf{x}-\int_{\Omega_\eps} W'(u_\eps)\varphi\, \de\textbf{x}=0,
\]
for all $\varphi\in H^1(\Omega_\eps)$.
We claim that
\begin{equation}\label{eq:EL_u_eps}
	\int_{\widetilde{\O}_\eps}\nabla v_\eps\cdot\nabla\psi\, \de\bfx=\delta^2\int_{\widetilde{\O}_\eps} W'(v_\eps)\psi\, \de\bfx,
\end{equation}
for every $\psi\in H^1(\widetilde{\O}_\eps)$.
Indeed, fix $\varphi\in H^1(\Omega_\eps)$. Using the change of variable $(\delta x',\delta y',\delta z')=(x,y,z)$, we get
\begin{align}\label{eq:change_gradient}
	\int_{\Omega_\eps}\nabla u_\eps \cdot \nabla\varphi\, \de\bold{x}
	=\delta^3\int_{\widetilde{\O}_\eps}\frac{1}{\delta^2}\nabla v_\eps\cdot\nabla\psi\, \de\bfx',
\end{align}
where $\psi(x,y,z)\coloneqq\varphi(\delta x,\delta y, \delta z)$.
Moreover,
\begin{align*}
	\int_{\Omega_\eps} W'(u_\eps)\varphi\, \de\bfx=\delta^3 \int_{\widetilde{\O}_\eps} W'(v_\eps)\psi\, \de\bfx'.
\end{align*}
This proves that $v_\eps$ is a weak solution to
\begin{align*}
	\begin{cases}
		\Delta v_\eps=\delta^2 W'(v_\eps)\quad&\text{in}\ \widetilde{\O}_\eps\,,\\[5pt]
		\displaystyle\frac{\partial v_\eps}{\partial \nu}=0\quad&\text{on}\ \partial\widetilde{\O}_\eps\,,
	\end{cases}
\end{align*}
as desired.

Now we want to obtain the limiting equation.
Since, by assumption, $(u_\eps)_\eps$ is uniformly bounded in $L^\infty$, the sequence $f_\eps\coloneqq \delta^2 W'(v_\eps)\chi_{\widetilde{\O}_\eps}$ converges strongly in $L^p_{\loc}(\mathbb{R}^3)$ to $f\coloneqq 0$, for all $p\geq 1$.
Moreover,
\[
\mathbb{R}^3\setminus \widetilde{\O}_\eps
	\stackrel{\mathrm{H}}{\longrightarrow} \mathbb{R}^3\setminus \Omega_\infty\,,\qquad\text{as $\eps\to0$.}
\]
Thus, using Theorem \ref{neumanconvergence}, we get that there exists $\hat{v}\in H^{1}_{\loc}(\O_\infty)$ such that, up to a subsequence,
\[
v_\eps\chi_{\widetilde{\O}_\eps}
	\to \hat{v}\chi_{\O_\infty} 
	\quad\text{ in } L^q_{\loc}(\mathbb{R}^3) \text{ for all } q\in[1,\infty),
\]
and
\[
\nabla v_\eps\chi_{\widetilde{\O}_\eps}
	\to \nabla \hat{v}\chi_{\O_\infty}
	\quad\text{ in } L^2_{\loc}(\mathbb{R}^3;\mathbb{R}^3),
\]
as $\eps\to0$.
This proves (ii).
In particular, if we fix $R>1$, we get that
\begin{equation}\label{eq:gradient_zero}
\int_{\Omega_\infty\cap B_{2R}}\abs{\nabla \hat{v}}^2\, \de\bold{x}
=\lim_{\eps\to0} \int_{\widetilde{\Omega}_\eps\cap B_{2R}}\abs{\nabla v_\eps}^2\, \de\bold{x} = 0,
\end{equation}
where we used \eqref{eq:EL_u_eps} with $v_\eps$ as a test function, together with the fact that $\| W'(v_\eps)v_\eps\|_{L^\infty}$ is uniformly bounded in $\eps$.

Therefore, recalling that $\O_\infty$ has two disjoint connected components, we get that there exist $m_1, m_2\in\mathbb{R}$ such that
\[
v_\eps \to m_1 \quad\text{ locally uniformly in } \O_\infty\cap B_{2R}\cap\{x<0\},
\]
and
\[
v_\eps \to m_2 \quad\text{ locally uniformly in } \O_\infty\cap B_{2R}\cap\{x>0\},
\]
as $\eps\to0$.
Moreover, for any choice of $a$ and $M$ such that $E(a,M)\subset B_{2R}$, we get that
\[
v_\eps \to m_1 \quad\text{ locally uniformly in } E^\ell(a,M)
\]
and
\[
v_\eps \to m_2 \quad\text{ locally uniformly in } E^\rr(a,M),
\]
as $\eps\to0$. 
Going back to the original coordinates gives us the desired result. We remark that, to obtain the locally uniform convergence, we may choose $p$ large enough in Theorem~\ref{neumanconvergence} so that the Sobolev embedding theorem can be invoked.

\noindent\emph{Step 1.2: boundary conditions on the external ellipsoid.}
Note that by assumption, we get that $u_\eps$ is a weak solution to
\begin{align*}
	\begin{cases}
		\Delta u_\eps= W'(u_\eps)\quad&\text{in}\ \O_\eps\,,\\[5pt]
		\displaystyle\frac{\partial u_\eps}{\partial \nu}=0\quad&\text{on}\ \partial\O_\eps\,.
	\end{cases}
\end{align*}
Moreover, $W'(u_\eps)\chi_{\O_\eps}\in L^p_{\loc}(\R{3})$ for all $p \geq 1$, since by assumption, $(u_\eps)_\eps$ is uniformly bounded in $L^\infty$.
Therefore, arguing as in the previous step, we get that $u_\eps\chi_{\O^\ell_\eps}$ converges uniformly to $\alpha$ and $u_\eps\chi_{\O^\rr_\eps}$ converges uniformly to $\beta$.
In particular, let $r_0>0$ be given by assumption (H3).
Define
\begin{align*}
E_{\eps,\rho}^\ell \coloneqq E_\rho\cap \{x<-\eps\}\qquad\text{and}\qquad
E_{\eps,\rho}^r\coloneqq  E_\rho\cap \{x>\eps\}.
\end{align*}
Then, from the above argument, we get that
\begin{align*}
	\nonumber\alpha-\gamma_\eps^\ell<&u_\eps<\alpha+\gamma_\eps^\ell\quad\text{on }\partial E_{\eps,\rho}^\ell\,, \\[5pt]
		\beta-\gamma_\eps^\rr<&u_\eps<\beta+\gamma_\eps^\rr\quad\text{on }\partial E_{\eps,\rho}^\rr\,,
\end{align*}
as $\eps\to0$, where
\begin{align*}
	\gamma_\eps^\ell\coloneqq ||u_\eps-\alpha||_{L^\infty(E_{\eps,\rho}^\ell )}\to0 \qquad\text{and}\qquad
	\gamma_\eps^r\coloneqq ||u_\eps-\beta||_{L^\infty( E_{\eps,\rho}^r )}\to0.
\end{align*}
This gives the desired result.\\

\blu{
\emph{Step 1.3: $m_1=m_2$}\,.
In order to prove that the two values are the same, we argue as follows.
First, we want to use the upper bound \eqref{upperboundueps}.
Note that both the optimal values for $A$ and $B$ given by \eqref{optimalAandB} converge to $(\alpha+\beta)/2$.
Therefore, from \eqref{upperboundueps} we get that
\begin{equation}\label{eq:upp_bound_finite}
\lim_{\eps\to0} \frac{|\ln(\eta/\delta)|}{\delta}F(u_\eps,N_\varepsilon) < +\infty.
\end{equation}
}
\blu{
Now, we want to get a lower bound as follow.
Without loss of generality, we can assume $m_1\leq m_2$.
Fix $\theta>0$, and let $\eps_0>0$ be given by Step 1.1.
Using \eqref{firstboundarycondition_1} and \eqref{firstboundarycondition_2}, for $\eps<\eps_0$, we estimate the energy inside the neck from below as follow:
\begin{align*}
F(u_\eps&, N_\eps) \geq
	\frac{1}{2}\int_{N_\eps}\abs{\nabla u_\eps}^2\, \de\bfx
	\\[5pt]
&\geq \inf\Big\{\frac{1}{2}\int_{N_\eps}\abs{\nabla v}^2\, \de\bfx
	: \ v\in H^1(N_\eps), v\leq  m_1+\theta  \text{ on } \partial N_\eps\cap\{x<0\}, \nonumber   \\[5pt]
&\hspace{4cm}
	\,\,
	v\geq m_2-\theta\text{ on } \partial N_\eps\cap\{x>0\} \Big\}  \nonumber \\[5pt]
&= \frac{\delta\eta}{\eps}(m_1-m_2+2\theta)^2.
\end{align*}
This lower bound yields that
\[
 \frac{|\ln(\eta/\delta)|}{\delta} F(u_\eps, N_\eps)
 	\geq |\ln(\eta/\delta)| \frac{\eta}{\eps} (m_1-m_2+2\theta)^2.
\]
Thus, taking the limit as $\theta\to0^+$ and recalling that in our regime
\[
\lim_{\eps\to0}|\ln(\eta/\delta)| \frac{\eta}{\eps} = +\infty,
\]
the only possibility in order to satisfy \eqref{eq:upp_bound_finite} is that $m_1=m_2$ as desired.\\
}

\emph{Step 1.4: lower bound of the energy}.
Fix $\gamma,\blu{\theta}>0$, and let $\varepsilon_0>0$ be the parameter given by Step 1.1 and Step 1.2.
\blu{Recall that Step 1.3 yields that $m_1=m_2$ in \eqref{firstboundarycondition_1} and \eqref{firstboundarycondition_2}. We will denote the common value by $m\in\R{}$.}
Then, for $\varepsilon<\varepsilon_0$, we have that
\begin{align}\label{eq:LR_eps}
F(u_\eps&, \O_\eps\setminus N_\eps) \geq
	\frac{1}{2}\int_{E^\ell(a_\eps, \rho)\setminus \overline{E}^\ell(a_\eps, M_\eps)}\abs{\nabla u_\eps}^2\, \de\bfx
	+ \frac{1}{2}\int_{E^\rr(a_\eps, \rho)\setminus \overline{E}^\rr(a_\eps, M_\eps)}\abs{\nabla u_\eps}^2\, \de\bfx \nonumber \\[5pt]
&\geq \inf\Big\{\frac{1}{2}\int_{E^\ell(a_\eps, \rho)\setminus \overline{E}^\ell(a_\eps, M_\eps)}\abs{\nabla v}^2\, \de\bfx
	: \ v\in H^1(E^\ell(a_\eps, \rho)\setminus \overline{E}^\ell(a_\eps, M_\eps)), \nonumber   \\[5pt]
&\hspace{4cm}
	v\leq m+\blu{\theta}\text{ on } \partial E^\ell(a_\eps, M_\eps),\,\,
	v\geq \alpha-\gamma\text{ on } \partial E^\ell(a_\eps, \rho) \Big\}  \nonumber \\[5pt]
&\hspace{0.5cm}+\inf\Big\{\frac{1}{2}\int_{E^\rr(a_\eps, \rho)\setminus \overline{E}^\rr(a_\eps, M_\eps)}\abs{\nabla v}^2\, \de\bfx
	:\ v\in H^1(E^\rr(a_\eps, \rho)\setminus \overline{E}^\rr(a_\eps, M_\eps)),  \nonumber  \\[5pt]
&\hspace{4cm}v \leq m+\blu{\theta}
		\text{ on } \partial E^\rr(a_\eps, M_\eps),\,\,
	v\geq \beta - \gamma \text{ on } \partial E^\rr(a_\eps, \rho)\Big\}  \nonumber  \\[5pt]
&= \inf\Big\{\frac{1}{2}\int_{E^\ell(a_\eps, \rho)\setminus \overline{E}^\ell(a_\eps, M_\eps)}\abs{\nabla v}^2\, \de\bfx
	: \ v\in H^1(E(a_\eps, \rho)^\ell\setminus \overline{E}^\ell(a_\eps, M_\eps)),  \nonumber  \\[5pt]
&\hspace{4cm}
	v= m+\blu{\theta} \text{ on } \partial E^\ell(a_\eps, M_\eps),\,\,
	v= \alpha-\gamma \text{ on } \partial E^\ell(a_\eps, \rho)
	\Big\}  \nonumber  \\[5pt]
&\hspace{0.5cm}+\inf\Big\{\frac{1}{2}\int_{E^\rr(a_\eps, \rho)\setminus \overline{E}^\rr(a_\eps, M_\eps)}\abs{\nabla v}^2\, \de\bfx
	:\ v\in H^1(E(a_\eps, \rho)^\rr\setminus \overline{E}^\rr(a_\eps, M_\eps)),  \nonumber  \\[5pt]
&\hspace{4cm}
	v = m+\blu{\theta}	\text{ on } \partial E^\rr(a_\eps, M_\eps),\,\,
	v= \beta - \gamma \text{ on } \partial E^\rr(a_\eps, \rho)\Big\}  \nonumber  \\[5pt]
&=: L_\eps + R_\eps.
\end{align}
\blu{Note that $L_\eps$ and $R_\eps$ depend also on $m, \gamma, \theta$, but we prefer to avoid specifying all of these parameters for the sake of notation.}

We now want to compute $L_\eps$ and $R_\eps$\,.
We show the argument for $L_\eps$\,. The result for $R_\eps$ will follow by using the same reasoning.
Arguing as in Section \ref{sec:elliptic_competitor}, we get that the solution to the problem defining $L_\eps$ is given, in prolate coordinates, by
\begin{align*}
w(\mu)=c \ln\abs{k \tanh(\mu/2)},
\end{align*}
where (see \eqref{eq:BS_elliptic_competitor}) \blu{
\begin{align*}
k=\frac{\exp\Big(\displaystyle\ln\abs{\frac{\tanh \rho}{\tanh M_\eps}}\frac{\blu{m+\theta}}{\alpha-\gamma-m-\blu{\theta}}\Big)}{\tanh \rho},\quad\quad\quad
c=\frac{\alpha-\gamma-m-\blu{\theta}}{\ln\Big(\displaystyle\frac{\tanh \rho}{\tanh M_\eps}\Big)}.
\end{align*}
}
In particular, (see \eqref{BoundA}), we get that
\begin{align}\label{loweboundwithgamma}
L_\eps = 	
	 \frac{\pi a_\eps(\blu{\alpha-\gamma-m-\theta})^2}{\ln\Big(\displaystyle\frac{\tanh \rho}{\tanh M_\eps}\Big)}.
\end{align}
Now, we want to understand the asymptotic behaviour of $L_\eps$.
First, we want to compute the asymptotic behaviour of the denominator on the right-hand side of \eqref{loweboundwithgamma}.
Note that the positions in \eqref{eq:a_e-M_e} yield that
\begin{equation}\label{eq:M_eps-a_eps}
	\tanh (2M_\eps) = \frac{\eta}{\delta}\qquad\text{and}\qquad
	a_\eps^2=4\delta^2-4\eta^2\approx 4\delta^2,\quad\text{as $\eps\to0$.}
\end{equation}
As consequence, we get the following asymptotic estimate
\begin{align}\label{lowerboundapprox}
\frac{1}{a_\eps}\ln\Big(\displaystyle\frac{\tanh \rho}{\tanh M_{\eps}}\Big)
=\frac{1}{a_\eps}\big( \ln \tanh \rho - \ln \tanh M_{\varepsilon}\big)
\approx\frac{\abs{\ln(\delta/\eta)}}{2\delta}.
\end{align}
Therefore, from \eqref{loweboundwithgamma}, and using the fact that we can take arbitrary $\gamma>0$ and $\delta>0$, we get
\begin{align*}
\liminf_{\eps\to0}\frac{\abs{\ln(\delta/\eta)}}{\delta}L_\eps
	\geq 2\pi  (m-\alpha)^2;
\end{align*}
in a similar way, we obtain that 
\begin{align*}
\liminf_{\eps\to0}\frac{\abs{\ln(\delta/\eta)}}{\delta}R_\eps
	\geq 2\pi (m-\beta)^2.
\end{align*}
\blu{These two estimates together imply \eqref{eq:window_lower_bound_energy}.}
\blu{
To conclude, we consider the function $f\colon \mathbb{R}\to\mathbb{R}$ given by
\begin{equation}\label{eq:f}
f(s)\coloneqq 2\pi[(s-\alpha)^2 + (s-\beta)^2].
\end{equation}
Then, the minimum of $f$ over the set $\alpha\leq s\leq \beta$ is given by
\[
f\left(\frac{\beta+\alpha}{2}\right)
	= \pi(\beta-\alpha)^2.
\]
Thus, from \eqref{eq:window_lower_bound_energy}, we obtain that
\[
\liminf_{\eps\to0} \frac{|\ln(\eta/\delta)|}{\delta}F(u_\eps,\O_\eps\setminus N_\eps)\geq \pi(\beta-\alpha)^2.
\]
This proves the claim.}\\

\noindent\emph{Step 2: energy estimates.}
\blu{Let $\xi_\eps$ be the function defined in \eqref{testfootball}.
Then, using \eqref{footballenergy} and \eqref{eq:window_lower_bound_energy}, we get that
\begin{align*}
\pi(\beta-\alpha)^2
	&=\lim_{\eps\to0} \frac{\left|\ln (\eta/\delta)\right|}{\delta} \frac{1}{2} \int_{E(a,M) \setminus E(a,m)} |\nabla \xi_\eps|^2 \, \de \bold{x}  \\
&\geq
	 \liminf_{\eps\to0} \frac{|\ln(\eta/\delta)|}{\delta}F(u_\eps,\O_\eps\setminus N_\eps) \geq \pi(\beta-\alpha)^2.
\end{align*}
This proves that
\[
\lim_{\eps\to0}\frac{\abs{\ln(\delta/\eta)}}{\delta}F(u_\eps,\O_\eps\setminus N_\eps)
	=\pi(\beta-\alpha)^2.
\]
In particular, since the minimizer of the function $f$ defined in \eqref{eq:f} is unique, this yields that
\[
m_1=m_2= \frac{\alpha+\beta}{2}.
\]
Therefore, using the result of Step 1.1, we get the validity of (i).}

\blu{
Moreover, by noticing that all the inequalities in \eqref{eq:LR_eps} are equalities, we get that
\[
\lim_{\eps\to0} \frac{|\ln(\eta/\delta)|}{\delta}F(u_\eps,\O_\eps)=
		\lim_{\eps\to0} \frac{|\ln(\eta/\delta)|}{\delta}F(u_\eps, \O_\eps\setminus N_\eps)=\pi(\beta-\alpha)^2.
\]
This proves the first part of the result.}\\

\noindent\emph{Step 3: limiting problem.}
\blu{Let $\hat{v}\in W^{2,p}(\Omega_\infty)$ be the function obtained in Step 1.2.}
First, we prove that $\hat{v}$ is an admissible competitor for the problem in (iii).
From (i), we know that
\[
\hat{v} = \frac{\alpha+\beta}{2},\qquad\text{ on } B_R\cap\Omega_\infty\,.
\]
We now prove that it satisfies also the boundary conditions at infinity.
Using the scale-invariant Poincar\'{e} inequality \eqref{poincare}, we get that
\[
\| v_\eps - \overline{v}_\eps \|_{L^6(\widetilde{\O}_\eps)} 
	\leq C \| \nabla v_\eps \|_{L^2(\widetilde{\O}_\eps)}\,,
\]
which, together with the fact that $\overline{v}_\eps\to \alpha\chi_{\O_\infty^\ell}+\beta\chi_{\O^\rr_\infty} \in L^6(\O_\infty)$, as $\eps\to0$, yields that $\hat{v}$ is an admissible competitor for the problem in (iii).

Finally, we prove that $\hat{v}$ solves the minimisation problem in (iii). The argument is similar to that of Step 3 of the proof of \cite[Theorem 4.1]{KohSla}.
Fix $M>|\alpha|,|\beta|$.
We can assume, without loss of generality, that every function $\varphi\in\mathcal{A}$ is such that $\|\varphi\|_{L^\infty(\O_\infty)}\leq M$. Indeed, given $\varphi\in\mathcal{A}$, by considering the truncation $\widetilde{\varphi}\coloneqq (\varphi \land M)\lor (-M)$ we get that $\widetilde{\varphi}\in\mathcal{A}$ and
\[
\frac{1}{2} \int_{\O_\infty} |\nabla \widetilde{\varphi}|^2\, \de \bold{x}
\leq
\frac{1}{2} \int_{\O_\infty} |\nabla \varphi|^2\,\de \bold{x}.
\]
Thus, let us take $\varphi\in\mathcal{A}$ with $\|\varphi\|_{L^\infty(\O_\infty)}\leq M$.
Define the function $\varphi_\eps:\Omega_\eps\to\R{}$ as
\[
\varphi_\eps(x,y,z)\coloneqq \varphi\left( \frac{x}{\delta}, \frac{y}{\delta}, \frac{z}{\delta} \right).
\]
Then, there exist constants $C,\widetilde{C}>0$, such that, for all $\eps>0$ it holds
\[
\| \varphi_\eps - u_{0,\eps} \|_{L^2(\O_\eps)}
\leq C \| \varphi_\eps - u_{0,\eps} \|_{L^6(\O_\eps)}
= \widetilde{C}\delta^{\frac{1}{2}}.
\]
Therefore, for $\eps$ sufficiently small, we get that $\varphi_\eps$ is an admissible competitor for the minimisation problem solved by $u_\eps$. Thus,
\begin{equation}\label{eq:min_F_v_eps}
F(u_\eps,\O_\eps)\leq F(\varphi_\eps,\O_\eps).
\end{equation}
Note that
\[
F(u_\eps,\O_\eps) = \frac{\delta}{2}\int_{\widetilde{\O}_\eps}
	|\nabla v_\eps|^2\ \de \textbf{x}
	+\int_{\O_\eps} W(u_\eps)\, \de \bold{x},
\]
and (recall that $\widetilde{\O}_\eps\subset\O_\infty$)
\[
F(\varphi_\eps,\O_\eps) = \frac{\delta}{2}\int_{\widetilde{\O}_\eps}
	|\nabla \varphi|^2\ \de \textbf{x}
	+\int_{\O_\eps} W(\varphi_\eps)\, \de \bold{x}.
\]
Thus, taking the liminf on both sides of \eqref{eq:min_F_v_eps}, and using the fact that $\varphi_\eps, u_\eps$ converges in $L^2$ to zeros of $W$, we get
\[
\frac{1}{2}\int_{\O_\infty}
	|\nabla \hat{v}|^2\, \de \bold{x}
\leq \liminf_{\eps\to0} F(u_\eps, \O_\eps)
\leq \liminf_{\eps\to0} F(\varphi_\eps,\O_\eps)
= \frac{1}{2}\int_{\O_\infty}
	|\nabla \varphi|^2\, \de \bold{x},
\]
\blu{where we used a change of variable to relate the Dirichlet energy of $u_\varepsilon$ with that of $\hat{v}$.}
This proves that that $\hat{v}$ solves the claimed minimisation problem.
This concludes the proof.
\end{proof}

\begin{remark}
We highlight that, from Step 1.1 of the proof (see \eqref{eq:gradient_zero}), it follows that the transition happens outside any ball of radius $\delta$ around the neck.
\end{remark}

\subsection{Narrow thick regime}
Here we consider the scaling of the energy in the narrow thick regime $\delta\gg\eps\approx\eta$, namely when
\begin{align*}
\lim_{\eps\to0}\frac{\eta}{\delta}=0\quad\text{and}\quad\lim_{\eps\to0} \frac{\eta}{\eps}=l\in(0,+\infty).
\end{align*}
Since the ellipsoidal competitor outside the neck provides an energy whose order is lower than the energy of the affine competitor in the neck, we expect the transition happening outside the neck. 
Denoting by $u_\eps$ a local minimiser of the functional \eqref{004}, the \blu{convenient rescaling} that allows us to both see a nice limiting space, and to use the rescaled Poincar\'{e} inequality is
\begin{align*}
	v_\eps(x,y,z)\coloneqq u_\eps(\delta x, \delta y, \delta z).
\end{align*}
Using this rescaling, the limiting domain becomes
\[
\widetilde{\Omega}_\infty = \{ x<0 \} \cup ( \{0\}\times [-1,1] \times \{0\} ) \cup \{ x>0 \}.
\]
However, note that, as $\eps\to0$,
\[
\R{3}\setminus \widetilde{\Omega}_\varepsilon
\stackrel{\mathrm{H}}{\rightarrow} 
\R{3}\setminus \Omega_\infty\,,
\]
where
\[
\widetilde{\Omega}_\varepsilon \coloneqq \frac{1}{\delta}\Omega_\varepsilon,\qquad
\Omega_\infty = \{ x<0 \} \cup \{ x>0 \}.
\]
The same argument used in the proof of Theorem \ref{Theoremwindowthick} yields the following result, therefore we omit the proof.

\begin{theorem}\label{Theoremnarrowthick}
Let $(u_\eps)_\eps$ be an admissible family of local minimisers as in Definition \ref{def:loc_min}.
Assume that $\delta\gg\eps\approx\eta$.
Then,
	\begin{equation*}
		\lim_{\eps\to0} \frac{|\ln(\eta/\delta)|}{\delta}F(u_\eps,\Omega_\eps)=
		\lim_{\eps\to0} \frac{|\ln(\eta/\delta)|}{\delta}F(u_\eps, \O_\eps\setminus N_\eps)=
		\pi(\beta-\alpha)^2. 
	\end{equation*}
Moreover, for $\varepsilon>0$ let $v_\varepsilon \colon \widetilde{\O}_\eps\to\R{}$ be defined as
\[
v_\eps(x,y,z)\coloneqq u_\eps(\delta x, \delta y, \delta z).
\]
Then, the following hold:
\begin{itemize}
\item[(i)] $v_\eps\to \frac{\alpha+\beta}{2}$ uniformly on $B_R$, as $\eps\to0$, \blu{for all $2R<r_0$, where $r_0$ is given by (H3)};
\item[(ii)] There exists $\hat{v}\in H^1_{\loc}(\Omega_\infty)$ such that $v_\eps\chi_{\widetilde{\O}_\eps}\to \hat{v}\chi_{\O_\infty}$ strongly in $H^1_{\loc}(\R{3})$, as $\eps\to0$;
\item[(iii)] The function $\hat{v}$ is the unique minimiser of the variational problem
\begin{align*}
\min\bigg\{ \frac{1}{2} \int_{\O_\infty} |\nabla v|^2\, \de x: v\in H^1_{\loc}(\O_\infty),\ &\, v-\alpha\chi_{\O_\infty^\ell}-\beta\chi_{\O^\rr_\infty} \in L^6(\O_\infty),\\
&\, v=\frac{\alpha+\beta}{2} \text{ on } B_R\cap\O_\infty  \bigg\},
	\end{align*}
where $\O^\ell_\infty\coloneqq \{x<0\}$, and $\O^\rr_\infty\coloneqq \{x>0\}$.
\end{itemize}
\end{theorem}

\subsection{Letter-box regime}
 We now consider the regime in which 
$ \delta \gg \eps \gg \eta$,
 namely when \begin{align*}
 	\lim_{\eps\to0} \frac{\delta}{\eps}=+\infty\quad\text{and}\quad\lim_{\eps\to0}\frac{\eta}{\delta}=0.
 \end{align*}
In this regime, the transition will happen all inside, all outside, or everywhere, depending on the parameter
\begin{align}\label{limitell}
	\ell\coloneqq \lim_{\eps\to 0}\frac{\delta\eta}{\eps}\frac{\ln\abs{\eta/\delta}}{\delta}.
\end{align}
In particular, we will prove that if $\ell\in(0,+\infty)$, then the transition happens everywhere, while if $\ell=0$ then transition occurs all inside and if  $\ell=+\infty$ all outside.

\subsubsection{Critical letter-box regime}
This sub-regime, corresponds to the case $\ell\in(0,+\infty)$.
We capture the transition in the bulk by applying a similar argument to the one in Theorem \ref{Theoremwindowthick}, in which around the neck, the rescaled profile $v_\eps$ converges to the average of the two phases, namely $(\alpha+\beta)/2$. In the critical letter-box regime we have instead that the rescaled profile converges to different constants $m_1$ on $\{x<0\}$ and $m_2$ on $\{x>0\}$. Once we have this information, we can understand how to describe the transition in the neck, by taking into account the fact that we know the boundary conditions. Then, a similar technique used in Theorem \ref{Theoremsuperthin} applies.

\begin{theorem}\label{Theoremcriticalletterbox}
Let $(u_\eps)_\eps$ be an admissible family of local minimisers as in Definition \ref{def:loc_min}.
Assume that $\delta \gg \eps \gg \eta$, and that
\[
\ell\coloneqq \lim_{\eps\to 0}\frac{\delta\eta}{\eps}\frac{\abs{\ln(\eta/\delta)}}{\delta}\in(0,+\infty).
\]
Then,
	\begin{align}\label{lbenergytotal}
		\lim_{\eps\to 0} \frac{\eps}{\delta\eta}F(u_\eps,\Omega_\eps)=\blu{\lim_{\varepsilon\to0}}\frac{\abs{\ln(\eta/\delta)}}{\delta\ell}F(u_\eps,\Omega_\eps)= \frac{\pi(\beta-\alpha)^2}{\pi+\ell}.
	\end{align}
In particular
\begin{align*}
\lim_{\eps\to 0} \frac{\eps}{\delta\eta}F(u_\eps,N_\eps)=\frac{\pi^2(\beta-\alpha)^2}{(\pi+\ell)^2},
\end{align*}
and
	\begin{equation*}
\lim_{\eps\to 0} \frac{\abs{\ln(\eta/\delta)}}{\delta}F(u_\eps,\O_\eps\setminus N_\eps)= \frac{\pi\ell^2(\beta-\alpha)^2}{(\pi+\ell)^2}.
\end{equation*}
Moreover:
\begin{enumerate}
	\item[$(i)$] Consider the rescaled profile $w_\eps:\overline{\O}_\eps\to \R{}$ defined as
	\begin{align}\label{letterboxneckrescaling}
		w_\eps(x,y,z)\coloneqq u_\eps(\eps x,\delta y, \eta z),
	\end{align} 
	where $\overline{\O}_\eps$ is the rescaled domain of $\O_\eps$. Let us assume that $w_\eps\wto \hat{w}$  in  $H^1(N)$, for some $\hat{w}\in H^1(N)$.
	Then, $\hat{w}(x,y,z)=w(x)$ where $w\in H^1([-1,1])$ is the unique minimiser of the variational problem
	\begin{align*}
		\min\bigg\{ \frac{1}{2} \int_{-1}^1 |v'|^2\,\de x:\ v\in H^1(-1,1),\ v(-1)&=\frac{\pi\alpha+\Big(\displaystyle\frac{\alpha+\beta}{2}\Big)\ell}{\pi+\ell},\\[5pt]
		v(1)&=\frac{\pi\beta+\Big(\displaystyle\frac{\alpha+\beta}{2}\Big)\ell}{\pi+\ell}\bigg\};
	\end{align*}
		\item[$(ii)$] Let 
		\begin{align*}
			\widetilde{\Omega}_\varepsilon \coloneqq \frac{1}{\delta}\Omega_\varepsilon,\quad\quad\quad\O_\infty\coloneqq \{x<0\}\cup \{x>0\}.
		\end{align*} 
		Define the rescaled profile $v_\eps\in H^1(\widetilde{\Omega}_\varepsilon)$ as
		\begin{align}\label{bulkrescalingletterbox}
			v_\eps(x,y,z)\coloneqq u_\eps(\delta x,\delta y, \delta z).
		\end{align}
		Then, there esists $\hat{v}\in H^1(\O_\infty)$ such that $v_\eps\chi_{\widetilde{\O}_\eps}\to \hat{v}\chi_{\O_\infty}$ strongly in $H^1_{\loc}(\R{3})$ as $\eps\to0$, where $\hat{v}$ is the solution of the minimisation problem
		\begin{align*}
		\min\bigg\{\frac{1}{2}\int_{\O_\infty}\abs{\nabla v}^2\, \de\bold{x}\ :\ v\in H^1(\O_\infty),\  &v-\alpha\chi_{\{x<0\}}-\beta\chi_{\{x>0\}}\in L^6(\O_\infty),\\[5pt]
		&v=\frac{\pi\alpha+\Big(\displaystyle\frac{\alpha+\beta}{2}\Big)\ell}{\pi+\ell}\text{ on } B_M\cap\{x<0\},\\[5pt]
		&v=\frac{\pi\beta+\Big(\displaystyle\frac{\alpha+\beta}{2}\Big)\ell}{\pi+\ell}\text{ on } B_M\cap\{x>0\}\bigg\},
		\end{align*}
		for some $M\geq 2$.
\end{enumerate}
\end{theorem}

\blu{
\begin{remark}
Notice that in part (i) we need to assume the weak convergence of the sequence $(w_\varepsilon)_\varepsilon$ to a limit function $\hat{w}$ in $H^1(N)$ because in this regime the energy bound is not strong enough to ensure compactness (see Remark \ref{rem:critical}).
\end{remark}
}
\begin{proof}[Proof of Theorem \ref{Theoremcriticalletterbox}]
First of all, note that, using \eqref{upperboundueps}, for any given constants $A,B\in\R{}$ with $A\leq B$, we have
\begin{align}\label{energyboundletterbox1}
	\frac{\abs{\ln(\eta/\delta)}}{\delta}	F(u_\varepsilon,\Omega_\eps)\leq \blu{2}\pi\big[(A-\alpha)^2+(B-\beta)^2\big]+\frac{\delta \eta}{\eps}\frac{\abs{\ln(\eta/\delta)}}{\delta}(B-A)^2.
\end{align}
Therefore, by \eqref{limitell}, for every $\lambda>0$ there exists $\eps_0>0$ such that for every $\eps<\eps_0$ we have
\begin{align}\label{energyboundletterbox2}
		\frac{\abs{\ln(\eta/\delta)}}{\delta}	F(u_\varepsilon,\Omega_\eps)\leq \blu{2}\pi\big[(A-\alpha)^2+(B-\beta)^2\big]+(\ell+\lambda)(B-A)^2.
\end{align}

\noindent\emph{Step 1: lower bound of the energy in the bulk.} The same strategy used in  Theorem \ref{Theoremwindowthick}, in which we obtained the boundary conditions at the edge of the neck, applies. \\
Consider the rescaling $v_\eps\in H^1(\widetilde{\Omega}_\eps)$ defined in \eqref{bulkrescalingletterbox} and the limiting domain $\O_\infty$.
By following the strategy in Step 1 of Theorem \ref{Theoremwindowthick} we obtain that there is $R> 0$ and $\hat{v}\in H^1(\O_\infty)$  such that $v_\eps\chi_{\widetilde{\O}_\eps}\to \hat{v}\chi_{\O_\infty}$ in $H^1_\loc(\R{3})$ and that $\hat{v}$ is constant on each connected component of $(\{x<0\}\cup\{x>0\})\cap B_R$. In other words, there are $m_1,m_2\in \R{}$ such that
\begin{align*}
	\hat{v}_{|B_R}=\begin{cases}
		m_1\quad\text{if }&x<0,\\[5pt]
		m_2\quad\text{if }&x>0.
	\end{cases}
\end{align*}
	Therefore, we can show that the following lower bound estimate holds (see \eqref{eq:window_lower_bound_energy})
	\begin{align}\label{letterbox:lbbulk}
	\liminf_{\eps\to0}\frac{\abs{\ln(\eta/\delta)}}{\delta}F(u_\eps,\O_\eps)&\geq\liminf_{\eps\to0}\frac{\abs{\ln(\eta/\delta)}}{\delta}F(u_\eps,\O_\eps\setminus N_\eps) \nonumber\\[1em]
	&\geq 2\pi\big[  (m_1-\alpha)^2+ (m_2-\beta)^2\big].
	\end{align}

\noindent\emph{Step 2: lower bound of the energy in the neck.} From the previous step, we infer that for any $\blu{\theta}>0$, there exists $\varepsilon_1>0$ such that, for every $\eps<\eps_1$,
	\begin{align*}
		m_1-\blu{\theta} \leq u_\eps \leq m_1+\blu{\theta}
		\quad\quad\text{on } \partial E^\ell(a_\eps, M_\eps),
	\end{align*}
	\begin{align*}
		m_2-\blu{\theta} \leq u_\eps \leq m_2+\blu{\theta}
		\quad\quad\text{on } \partial E^\rr(a_\eps, M_\eps),
	\end{align*}
	for a suitable ellipsoid $E(a_\eps,M_\eps)$. From that, we obtain
	 a lower bound of the energy in the neck. Indeed,
	 \begin{align}\label{letterbox:lbneck}
	 	\nonumber\frac{\eps}{2\delta\eta}\int_{N_\eps}\abs{\nabla u_\eps}^2\, \de\bfx &\geq\inf\Big\{\frac{1}{2} \int_{N_\eps}\abs{\nabla v}^2\, \de\bfx:\ v\in H^1(N_\eps),\\[5pt] \nonumber&\hspace{1.5cm}v\geq m_1-\blu{\theta} \text{ on } \{x=-\eps\}
	 	\text{ and }v\leq m_2+\blu{\theta} \text{ on } \{x=\eps\}\Big\}\\[5pt]
	 	\nonumber&\geq\inf\Big\{\frac{1}{2} \int_{N_\eps}\abs{\nabla v}^2\, \de\bfx:\ v\in H^1(N_\eps),\\[5pt] \nonumber&\hspace{1.5cm}v=m_1-\blu{\theta} \text{ on } \{x=-\eps\}
	 	\text{ and }v=m_2+\blu{\theta} \text{ on } \{x=\eps\}\Big\}\\[5pt]
	 	&=\big(m_1-m_2-2\blu{\theta}\big)^2,
	 \end{align}
	 where in the last step we used the fact that the minimiser of the above minimisation problem is given by the affine function.
	 
\noindent\emph{Step 3: limit of the energy.} By putting together \eqref{letterbox:lbbulk}, \eqref{letterbox:lbneck} and making use of \eqref{limitell}, we obtain
	\begin{align}\label{splittingenergyletterbox}
	\nonumber	\liminf_{\eps\to0}\frac{\abs{\ln(\eta/\delta)}}{2\delta}\int_{\O_\eps}\abs{\nabla u_\eps}^2\,\de\bold{x}&\geq\liminf_{\eps\to0}\frac{\abs{\ln(\eta/\delta)}}{2\delta}\int_{\O_\eps^\ell\cup\O_\eps^\rr}
	\abs{\nabla u_\eps}^2\de\textbf{x}\\[5pt]
		\nonumber&\hspace{1cm}+\liminf_{\eps\to0}\frac{\abs{\ln(\eta/\delta)}}{2\delta}\int_{N_\eps
		}\abs{\nabla u_\eps}^2\,\de\bold{x}\\[5pt] 
		&\geq 2\pi\big[  (m_1-\alpha)^2+ (m_2-\beta)^2\big]+\ell(m_1-m_2-2\blu{\theta})^2.
	\end{align}
On the other hand, from \eqref{energyboundletterbox2}, we have
\begin{align*}
	2\pi\big[  (m_1-\alpha)^2+ (m_2-\beta)^2\big]+(\ell+\lambda)\big(m_1-m_2\big)^2
	\geq\limsup_{\eps\to0}\frac{\abs{\ln(\eta/\delta)}}{2\delta}\int_{\O_\eps}\abs{\nabla u_\eps}^2\,\de\bold{x}.
\end{align*}
By letting $\blu{\theta},\lambda\to0$ in the above two inequalities, we obtain
	\begin{align}\label{letterboxtooptimize}
		\lim_{\eps\to0}\frac{\abs{\ln(\eta/\delta)}}{2\delta}\int_{\O_\eps}\abs{\nabla u_\eps}^2\,\de\bold{x}
		=
		2\pi\big[  (m_1-\alpha)^2+ (m_2-\beta)^2\big]+\ell(m_1-m_2)^2.
	\end{align}
The right-hand side is minimized for
	\begin{align}\label{lbm1m2}
m_1=\frac{\pi\alpha+\Big(\displaystyle\frac{\alpha+\beta}{2}\Big)\ell}{\pi+\ell}\quad\text{and}\quad m_2=\frac{\pi\beta+\Big(\displaystyle\frac{\alpha+\beta}{2}\Big)\ell}{\pi+\ell},
	\end{align}
which gives 
	\begin{align*}
	\lim_{\eps\to 0} \frac{\eps}{\delta\eta}F(u_\eps,\O_\eps)= \frac{\pi(\beta-\alpha)^2}{\pi+\ell}.
\end{align*}
In particular, by noticing that all the inequalities in \eqref{letterbox:lbbulk}, \eqref{letterbox:lbneck}, and \eqref{splittingenergyletterbox} are equalities, we get
	\begin{align*}
		\lim_{\eps\to 0} \frac{\eps}{\delta\eta}F(u_\eps,N_\eps)= \frac{\pi^2(\beta-\alpha)^2}{(\pi+\ell)^2},
	\end{align*}
	and
	\begin{align*}
		\lim_{\eps\to 0} \frac{\abs{\ln(\eta/\delta)}}{\delta}F(u_\eps,\O_\eps\setminus N_\eps)= \frac{\pi\ell^2(\beta-\alpha)^2}{(\pi+\ell)^2}.
	\end{align*}

\noindent\emph{Step 4: limiting problems.} Now we investigate the variational problem that the rescalings $v_\eps$, defined in \eqref{bulkrescalingletterbox}, in the bulk  and $w_\eps$, defined in \eqref{letterboxneckrescaling}, in the neck satisfy asymptotically.

\noindent\emph{Step 4.1: limiting problem in the neck.} By acting like Step 1.1 of Theorem \ref{Theoremwindowthick}, let $R>0$ be such that, as $\eps\to0$,
\begin{align*}
	u_\eps&\to m_1\quad\text{uniformly on } B_{\delta R}\cap \{x\leq-\eps\},\\[5pt]
		u_\eps&\to m_2\quad\text{uniformly on } B_{\delta R}\cap \{x\geq\eps\},
\end{align*} 
with $m_1,m_2$ defined in \eqref{lbm1m2}.
In particular, $u_\eps$ has asymptotic boundary conditions ad the edge of the neck $m_1$ and $m_2$ respectively.
Using the fact that $\delta\gg \eps\gg \eta$, it follows that, if we consider the rescaling \eqref{letterboxneckrescaling},
\begin{equation}\label{uniformconvergenceofweps}
	\begin{aligned}
		w_\eps &\to m_1 \quad \text{uniformly on } \{x=-1\}\times[-1,1]^2, \\[5pt]
		w_\eps &\to m_2 \quad \text{uniformly on } \{x=1\}\times[-1,1]^2,
	\end{aligned}
\end{equation}
as $\eps\to0$, which gives us the asymptotic boundary conditions at the edge of the neck satisfied by the limiting profile.
By assumption, there exists $\hat{w}\in H^1(N)$ such that $w_\eps\wto \hat{w}$ in $H^1(N)$ as $\eps\to0$ and, from \eqref{uniformconvergenceofweps}, $\hat{w}$ is an admissible competitor for the variational problem in $(i)$. Moreover,
\begin{align}\label{letterboxw}
\nonumber(m_1-m_2)^2&\geq \frac{\eps}{2\delta\eta}\liminf_{\eps\to0}\int_{N_\eps} \abs{\nabla u_\eps}^2\, \de\bfx \\[5pt]
&= \liminf_{\eps\to0}\frac{1}{2} \int_N \Big((\partial_x w_\eps)^2+ \frac{\eps^2}{\delta^2}(\partial_y w_\eps)^2+ \frac{\eps^2}{\eta^2}(\partial_z w_\eps)^2\Big)\, \de\bfx.
\end{align} 
Since $\eps/\eta\to\infty$, we have that $\hat{w}$ do not depend on the variable $z$ and since we know that $w_\eps\wto\hat{w}$ we have, from \eqref{letterboxw}, that
\begin{align*}
(m_1-m_2)^2&\geq\liminf_{\eps\to0}\frac{1}{2}
	\int_N (\partial_x w_\eps)^2\,\de\bfx 
	=\int_{[-1,1]^2} (\partial_x \hat{w}(x,y))^2\,\de x \de y\\[5pt]
&=\int_{-1}^1 \int_{-1}^1 (\partial_x \hat{w}(x,y))^2\,\de x \de y
\geq \frac{1}{2}\int_{-1}^1 \left|\int_{-1}^1 \partial_x \hat{w}(x,y)\,\de x \right|^2 \de y\\[5pt]
&=\frac{1}{2}\int_{-1}^1 \big|\hat{w}(1,y)-\hat{w}(-1,y)\big|^2 \de y\\[5pt]
&=(m_1-m_2)^2,
\end{align*}
where in the second inequality we used Jensen inequality.
Therefore, we conclude that $\hat{w}(x,y,z)=w(x)$ for $w\in H^1([-1,1])$ and $w$ solves the variational problem in $(i)$.

\noindent\emph{Step 4.2: limiting problem.} We use a similar argument to the one employed in Step 3 of Theorem \ref{Theoremwindowthick} applies. More specifically, $\hat{v}$ is admissible competitor for the problem in $(ii)$ and
\begin{align*}
	\hat{v}_{|B_R\cap\O_\infty}=\begin{cases}
\displaystyle\frac{\pi\alpha+\Big(\displaystyle\frac{\alpha+\beta}{2}\Big)\ell}{\pi+\ell}\quad\text{on } &\{x<0\},\\[10pt] \displaystyle\frac{\pi\beta+\Big(\displaystyle\frac{\alpha+\beta}{2}\Big)\ell}{\pi+\ell}\quad\text{on } &\{x>0\}.
	\end{cases}
\end{align*} 
By using the rescaled Poincar\'{e} inequality, we get that
\[
\| v_\eps - \overline{v}_\eps \|_{L^6(\widetilde{\O}_\eps)} 
\leq C \| \nabla v_\eps \|_{L^2(\widetilde{\O}_\eps)}.
\]
In analogy with Step 2 of Theorem \ref{Theoremsuperthin} we get $\overline{v}_\eps\to \alpha\chi_{\O_\infty^\ell}+\beta\chi_{\O^\rr_\infty} \in L^6(\O_\infty)$, as $\eps\to0$.
Finally, by applying the last part of Step 3 of Theorem \ref{Theoremwindowthick}, we have that $\hat{v}$ solves the variational problem in $(ii)$.
 \end{proof}

\begin{remark}
Note that in this case, we do not have compactness of the rescaled profile $w_\eps$ inside the neck. This is due to the fact that the \blu{chosen} rescaling, that allows us to see the neck at scale one, does not give a uniform bound on the gradient of the rescaled profile (in particular, the derivative with respect to the variable $y$ cannot be bounded).
\end{remark}

We now investigate the remaining two sub-regimes.

\subsubsection{Super-critical Letter-box regime}
In this sub-regime, we have
\begin{align}
\label{linfty}
\ell= \lim_{\eps\to 0}\frac{\delta\eta}{\eps}\frac{\abs{\ln(\eta/\delta)}}{\delta}=+\infty.
\end{align}
In this case, we recover the same result as in Theorem \ref{Theoremwindowthick}.
\begin{theorem}\label{Theoremsupercritical}
	Let $(u_\eps)_\eps$ be an admissible family of local minimisers as in Definition \ref{def:loc_min}.
	Assume that $\delta \gg \eps \gg \eta$ and that $\ell=+\infty$.
	Then,
	\begin{equation*}
		\lim_{\eps\to0} \frac{|\ln(\eta/\delta)|}{\delta}F(u_\eps,\O_\eps)=
		\lim_{\eps\to0} \frac{|\ln(\eta/\delta)|}{\delta}F(u_\eps, \O_\eps\setminus N_\eps)=
		\pi(\beta-\alpha)^2.
	\end{equation*}
	Moreover, let 
	\begin{align*}
		\widetilde{\Omega}_\varepsilon \coloneqq \frac{1}{\delta}\Omega_\varepsilon,\quad\quad\quad\O_\infty\coloneqq \{x<0\}\cup \{x>0\}.
	\end{align*}
	 Define the rescaled profile $v_\varepsilon \colon \widetilde{\O}_\eps\to\R{}$ be defined as
	\[
	v_\eps(x,y,z)\coloneqq u_\eps(\delta x, \delta y, \delta z).
	\]
	Then, that the following hold:
	\begin{itemize}
		\item[(i)] $v_\eps\to \frac{\alpha+\beta}{2}$ uniformly on $B_R$, as $\eps\to0$, \blu{for all $2R<r_0$, where $r_0$ is given by (H3)};
		\item[(ii)] There exists $\hat{v}\in H^1_{\loc}(\Omega_\infty)$ such that $v_\eps\chi_{\widetilde{\O}_\eps}\to \hat{v}\chi_{\O_\infty}$ strongly in $H^1_{\loc}(\R{3})$, as $\eps\to0$;
		\item[(iii)] The function $\hat{v}$ is the unique minimizer of the variational problem
		\begin{align*}
			\min\bigg\{ \frac{1}{2} \int_{\O_\infty} |\nabla v|^2\, \de x:
			v\in \mathcal{A}  \bigg\},
		\end{align*}
		where
		\[
		\mathcal{A}\coloneqq \Bigl\{v\in H^1_{\loc}(\O_\infty),\ v-\alpha\chi_{\O_\infty^\ell}-\beta\chi_{\O^\rr_\infty} \in L^6(\O_\infty),\, v=\frac{\alpha+\beta}{2} \text{ on } B_R\cap\O_\infty \Bigr\},
		\]
		and $\O^\ell_\infty\coloneqq \{x<0\}$, and $\O^\rr_\infty\coloneqq \{x>0\}$.
	\end{itemize}
\end{theorem}
\begin{proof}
The proof follows the ones for Theorems \ref{Theoremcriticalletterbox} and \ref{Theoremwindowthick}.

\noindent\emph{Step 1: bound of the energy.} By \eqref{upperboundueps}, for any given constants $A,B\in\mathbb{R}$ with $\alpha\leq A\leq B\leq\beta$, we have
\begin{align*}
	\frac{\abs{\ln(\eta/\delta)}}{\delta}F(u_\eps,\Omega_\eps)\leq \blu{2} \pi\big[(A-\alpha)^2+(B-\beta)^2 \big]
	+\frac{\delta \eta}{\eps}\frac{\abs{\ln(\eta/\delta)}}{\delta}(B-A)^2.
\end{align*}
Since
\[
\lim_{\eps\to0} \frac{\delta \eta}{\eps}\frac{\abs{\ln(\eta/\delta)}}{\delta} =+\infty,
\]
the only way to get that the right-hand side of the above inequality is bounded uniformly in $\eps$, is to choose $A=B$. This gives the estimate
\begin{align*}
	\frac{\abs{\ln(\eta/\delta)}}{\delta}F(u_\eps,\Omega_\eps)\leq \blu{2} \pi[(A-\alpha)^2+(A-\beta)^2].
\end{align*}

\noindent\emph{Step 2. Lower bound of the energy.} From \eqref{letterboxtooptimize}, we get
	\begin{align*}
	\nonumber	\liminf_{\eps\to0}\frac{\abs{\ln(\eta/\delta)}}{2\delta}\int_{\O_\eps}\abs{\nabla u_\eps}^2\,\de\bold{x}&\geq\liminf_{\eps\to0}\frac{\abs{\ln(\eta/\delta)}}{2\delta}\int_{\O_\eps^\ell\cup\O_\eps^\rr}
	\abs{\nabla u_\eps}^2\,\de\bold{x}\\[5pt]
	\nonumber&\hspace{1cm}+\liminf_{\eps\to0}\frac{\abs{\ln(\eta/\delta)}}{2\delta}\int_{N_\eps
	}\abs{\nabla u_\eps}^2\,\de\bold{x}\\[5pt] 
	\nonumber&\geq\liminf_{\eps\to0}\Big[2\pi\big[  (m_1-\alpha)^2+ (m_2-\beta)^2\big] \\
	&\phantom{\liminf_{\eps\to0}\Big[}\quad+\frac{\abs{\ln(\eta/\delta)}}{\delta}\frac{\delta\eta}{\eps}(m_1-m_2)^2\Big].
\end{align*}
Therefore, as in Theorem \ref{Theoremwindowthick}, we have an optimality condition on $m_1$ and $m_2$, which, together with \eqref{linfty}, leads to
\[
	A=m_1= m_2=\frac{\alpha+\beta}{2}.
\]
Therefore
\begin{align*}
\lim_{\eps\to0} \frac{\abs{\ln(\eta/\delta)}}{\delta}F(u_\eps,\O_\eps)=
\lim_{\eps\to0} \frac{|\ln(\eta/\delta)|}{\delta}F(u_\eps, \O_\eps\setminus N_\eps)=
\pi(\beta-\alpha)^2.
\end{align*}
The rest of the proof is identical to the one in Theorems \ref{Theoremwindowthick} and \ref{Theoremcriticalletterbox} and we obtain the desired result.
\end{proof}

\subsubsection{Sub-critical Letter-box regime}
In this sub-regime, we have
\begin{align}
\label{ellsubcritical}
\ell= \lim_{\eps\to 0}\frac{\delta\eta}{\eps}\frac{\abs{\ln(\eta/\delta)}}{\delta}=0.
\end{align}
In this case, we recover the same result as in Theorem \ref{Theoremsuperthin}.

\begin{theorem}
	Let $(u_\eps)_\eps$ be an admissible family of local minimisers as in Definition \ref{def:loc_min}.
		Assume that $\delta \gg \eps \gg \eta$ and that $\ell=0$.
	Then,
	\begin{align*}
		\lim_{\eps\to0} \frac{\eps}{\delta\eta}F(u_\eps,\O_\eps)
		=\lim_{\eps\to0} \frac{\eps}{\delta\eta}F(u_\eps, N_\eps)
		= (\beta-\alpha)^2.
	\end{align*}
Define the rescaled profile $v_\varepsilon \colon \overline{\O}_\eps\to\R{}$ as
	\[
	v_\eps(x,y,z)\coloneqq u_\eps(\eps x,\delta y, \eta z),
	\]
	where $\overline{\O}_\eps$ is the rescaled domain of $\O_\eps$.
	Then the following hold:
	\begin{itemize}
		\item[(i)] $v_\eps\chi_{\widetilde{\O}^\ell_\eps}- \alpha\chi_{\widetilde{\O}^\ell_\eps}\to0$ in $L^6(\R{3})$ and
		$v_\eps\chi_{\widetilde{\O}^\rr_\eps}- \beta\chi_{\widetilde{\O}^\rr_\eps}\to0$ in $L^6(\R{3})$, as $\eps\to0$;
		\item[(ii)] There exists $\hat{v}\in H^1(N)$ such that $v_\eps\wto \hat{v}$ weakly in $H^1(N)$, as $\eps\to0$; 
		\item[(iii)] It holds that $\hat{v}(x,y,z)=v(x)$, where $v\in H^1(-1,1)$ is the unique minimizer of the variational problem
		\begin{align*}
			\min\bigg\{ \frac{1}{2} \int_{-1}^1 |v'|^2\, \de x:\ v\in H^1(-1,1),\ v(-1)=\alpha,\ v(1)=\beta\bigg\},
		\end{align*}
		In particular, $v(x)=\dfrac{\beta-\alpha}{2}x+\dfrac{\alpha+\beta}{2}$.
\end{itemize}
\end{theorem}

\begin{proof}
	The proof is an adaptation of Theorem \ref{Theoremsupercritical} and we remark only the differences. Regarding the bound of the energy, we have that for any given constants $A,B>0$ with $A\leq B$, we have
	\begin{align*}
		\frac{\eps}{\delta\eta}F(u_\eps,\Omega_\eps)\leq \blu{2} \pi\frac{\eps}{\delta\eta}\frac{\delta}{\abs{\ln(\eta/\delta)}} \big[(A-\alpha)^2+(B-\beta)^2\big]+(B-A)^2.
	\end{align*} 
Since \eqref{ellsubcritical} holds, we obtain a bound for the energy which is compatible with a transition inside the neck by choosing $A=\alpha$ and $B=\beta$.

By using the same techniques as in Theorems \eqref{Theoremwindowthick} and \eqref{Theoremcriticalletterbox}, we obtain the desired result.
\end{proof}

\section*{Acknowledgments}
This research was carried out at the departments of Mathematics of Politecnico di Torino and Radboud University, whose hospitality is gratefully acknowledged. RC and MM are members of GNAMPA (INdAM).
MM thanks Gurgen Hayrapetyan for useful discussions.

\label{page:e}

\begin{thebibliography}{99}
\bibitem{ArrCar04}
{\sc J.~M. Arrieta and A.~N. Carvalho}, {\em Spectral convergence and nonlinear
  dynamics of reaction-diffusion equations under perturbations of the domain},
  J. Differential Equations, 199 (2004), pp.~143--178.

\bibitem{ArrCarLC06}
{\sc J.~M. Arrieta, A.~N. Carvalho, and G.~Lozada-Cruz}, {\em Dynamics in
  dumbbell domains. {I}. {C}ontinuity of the set of equilibria}, J.
  Differential Equations, 231 (2006), pp.~551--597.

\bibitem{ArrCarLC09a}
\leavevmode\vrule height 2pt depth -1.6pt width 23pt, {\em Dynamics in dumbbell
  domains. {II}. {T}he limiting problem}, J. Differential Equations, 247
  (2009), pp.~174--202.

\bibitem{ArrCarLC09b}
\leavevmode\vrule height 2pt depth -1.6pt width 23pt, {\em Dynamics in dumbbell
  domains. {III}. {C}ontinuity of attractors}, J. Differential Equations, 247
  (2009), pp.~225--259.

\bibitem{Bru99}
{\sc P.~Bruno}, {\em Geometrically constrained magnetic wall}, Phys. Rev.
  Lett., 83 (1999), pp.~2425--2428.

\bibitem{CabFreMorPer01}
{\sc E.~Cabib, L.~Freddi, A.~Morassi, and D.~Percivale}, {\em Thin notched
  beams}, J. Elasticity, 64 (2001), pp.~157--178.

\bibitem{CDLLMur04b}
{\sc J.~Casado-D\'iaz, M.~Luna-Laynez, and F.~c. Murat}, {\em Asymptotic
  behavior of an elastic beam fixed on a small part of one of its extremities},
  C. R. Math. Acad. Sci. Paris, 338 (2004), pp.~975--980.

\bibitem{CDLLMur04}
\leavevmode\vrule height 2pt depth -1.6pt width 23pt, {\em Asymptotic behavior
  of diffusion problems in a domain made of two cylinders of different
  diameters and lengths}, C. R. Math. Acad. Sci. Paris, 338 (2004),
  pp.~133--138.

\bibitem{CDLLMur08}
\leavevmode\vrule height 2pt depth -1.6pt width 23pt, {\em The diffusion
  equation in a notched beam}, Calc. Var. Partial Differential Equations, 31
  (2008), pp.~297--323.

\bibitem{ChaDov}
{\sc A.~Chambolle and F.~Doveri}, {\em Continuity of neumann linear elliptic
  problems on varying two-dimensional bounded open sets}, Commun. Partial
  Differ. Equ., 22 (1997), pp.~811--840.

\bibitem{CheYan10}
{\sc S.~Chen and Y.~Yang}, {\em Phase transition solutions in geometrically
  constrained magnetic domain wall models}, Journal of Mathematical Physics, 51
  (2010), p.~023504.

\bibitem{DMEboPon}
{\sc G.~Dal~Maso, F.~Ebobisse, and M.~Ponsiglione}, {\em A stability result for
  nonlinear {N}eumann problems under boundary variations}, J. Math. Pures
  Appl., 82 (2003), pp.~503--532.

\bibitem{HalVeg84}
{\sc J.~K. Hale and J.~Vegas}, {\em A nonlinear parabolic equation with varying
  domain}, Arch. Rational Mech. Anal., 86 (1984), pp.~99--123.

\bibitem{Jim84}
{\sc S.~Jimbo}, {\em Singular perturbation of domains and semilinear elliptic
  equation}, J. Fac. Sci. Univ. Tokyo Sect. IA Math., 35 (1988), pp.~27--76.

\bibitem{Jim88}
\leavevmode\vrule height 2pt depth -1.6pt width 23pt, {\em Singular
  perturbation of domains and the semilinear elliptic equation. {II}}, J.
  Differential Equations, 75 (1988), pp.~264--289.

\bibitem{Jim04}
\leavevmode\vrule height 2pt depth -1.6pt width 23pt, {\em Singular
  perturbation of domains and semilinear elliptic equations. {III}}, Hokkaido
  Math. J., 33 (2004), pp.~11--45.

\bibitem{JubAllBis04}
{\sc P.-O. Jubert, R.~Allenspach, and A.~Bischof}, {\em Magnetic domain walls
  in constrained geometries}, Phys. Rev. B, 69 (2004), p.~220410.

\bibitem{Klaui_2008}
{\sc M.~Kläui}, {\em Head-to-head domain walls in magnetic nanostructures},
  Journal of Physics: Condensed Matter, 20 (2008), p.~313001.

\bibitem{KohSla}
{\sc R.~V. Kohn and V.~V. Slastikov}, {\em Geometrically constrained walls},
  Calc. Var. Partial Differential Equations, 28 (2007), pp.~33--57.

\bibitem{Mol_etal02}
{\sc V.~A. Molyneux, V.~V. Osipov, and E.~V. Ponizovskaya}, {\em Stable two-
  and three-dimensional geometrically constrained magnetic structures: The
  action of magnetic fields}, Phys. Rev. B, 65 (2002), p.~184425.

\bibitem{MolOsiPon02}
\leavevmode\vrule height 2pt depth -1.6pt width 23pt, {\em Stable two- and
  three-dimensional geometrically constrained magnetic structures: The action
  of magnetic fields}, Phys. Rev. B, 65 (2002), p.~184425.

\bibitem{MorSla12}
{\sc M.~Morini and V.~Slastikov}, {\em Geometrically constrained walls in two
  dimensions}, Arch. Ration. Mech. Anal., 203 (2012), pp.~621--692.

\bibitem{MorSla15}
\leavevmode\vrule height 2pt depth -1.6pt width 23pt, {\em Geometrically
  induced phase transitions in two-dimensional dumbbell-shaped domains}, J.
  Differential Equations, 259 (2015), pp.~1560--1605.

\bibitem{RubSchSte04}
{\sc J.~Rubinstein, M.~Schatzman, and P.~Sternberg}, {\em Ginzburg-{L}andau
  model in thin loops with narrow constrictions}, SIAM J. Appl. Math., 64
  (2004), pp.~2186--2204.

\bibitem{Sasaki_etal12}
{\sc M.~Sasaki, K.~Matsushita, J.~Sato, and H.~Imamura}, {\em Thermal stability
  of the geometrically constrained magnetic wall and its effect on a
  domain-wall spin valve}, Journal of Applied Physics, 111 (2012), p.~083903.

\bibitem{Tatara_etal99}
{\sc G.~Tatara, Y.-W. Zhao, M.~Mu\~noz, and N.~Garc\'{\i}a}, {\em Domain wall
  scattering explains 300\% ballistic magnetoconductance of nanocontacts},
  Phys. Rev. Lett., 83 (1999), pp.~2030--2033.
\end{thebibliography}
\end{document}